\documentclass[10pt,reqno]{article}
\usepackage{amsmath, amsthm, amssymb, stmaryrd}
\usepackage{hyperref}
\usepackage{enumerate}
\usepackage{url}
\usepackage{color}
\usepackage{tikz}
\usetikzlibrary{patterns}
\usepackage[T1,T2A]{fontenc}
\usepackage[utf8]{inputenc}

\newtheorem{maintheorem}{Theorem}
\newtheorem{thm}{Theorem}
\newtheorem{lemma}[thm]{Lemma}

\theoremstyle{definition}

\theoremstyle{remark}
\newtheorem{remark}{Remark}

\def \le {\leqslant}
\def \ge {\geqslant}

\usepackage{tikz} 
\usetikzlibrary{calc,arrows,decorations.pathreplacing,fadings,3d,positioning}

\begin{document}
\title{On geometry of simultaneous approximation to three real numbers}
\author{by {\bf Antoine Marnat } and {\bf Nikolay Moshchevitin}}

\date{}
\maketitle

{\bf Abstract:} Considering simultaneous approximation to three numbers, we study the geometry of the sequence of best approximations. We provide a sharper lower bound for the ratio between ordinary and uniform exponent of Diophantine approximation, optimal in terms of this geometry.

 {\bf AMS Subject Classification:} 11J13
 
 {\bf Keywords:} simultaneous approximation, Diophantine exponents, best approximations. 

\section{ Simultaneous approximation and Diophantine exponents}

In this paper we deal with simultaneous approximation to real numbers $\alpha_1,\ldots,\alpha_n$,
that is integer solutions  $(q,a_1,\ldots,a_n)\in \mathbb{Z}^{n+1}$ of the system  inequalities
\begin{equation}\label{s4}
\begin{cases}
\displaystyle{\max_{1\le j \le n}} \, |q\alpha_j - a_j| \le Q^{-\gamma},
\cr
1\le q \le Q.
\end{cases}
\end{equation}
The \emph{ordinary  Diophantine exponent}  $ \omega = \omega(\pmb{\alpha})$ for real vector 
$\pmb{\alpha}\in \mathbb{R}^n$ is defined as supremum  over all $\gamma$ such that the system \eqref{s4} has integer solutions for arbitrary large values of $Q$.
The {\it uniform   Diophantine exponent}  $ \hat{\omega} = \hat{\omega}(\pmb{\alpha}) $ is defined as supremum  over all $\gamma$ such that the system \eqref{s4} has integer solutions for all $Q$ large enough.
It is well known that for $\pmb{\alpha}\not\in \mathbb{Q}^n$ the value of  $ \hat{\omega} =  \hat{\omega}(\pmb{\alpha}) $ always satisfies the inequalities
$$
\frac{1}{n} \le \hat{\omega}
\le 1.
$$
We need to distinguish  cases $ \hat{\omega}<1 $ and $\hat{\omega} = 1$.
Generalizing  the results by  Jarn\'{\i}k \cite{J}, in \cite{MaMo}
it was proven that under the condition that 
the numbers $1,\alpha_1,\ldots,\alpha_n$ are linearly independent over $\mathbb{Q}$ one has 
\begin{equation}\label{MaMo}
\frac{\omega}{\hat{\omega}} \ge G_n(\hat{\omega}),
\end{equation}
where $G_n(\hat{\omega})$
is the unique positive root of the equation
$$
g^{n-1} = \frac{\hat{\omega}}{1-\hat{\omega}} (g^{n-2}+ \cdots +g+1)
$$
in the case $ \hat{\omega} < 1$ and 
$G_n(\hat{\omega})=\infty $
in the case $ \hat{\omega} = 1$.

An alternative proof of this result is given in  \cite{Nguen}. This bound is known to be optimal.\\

In this paper we deal with approximations to three numbers. In this cas, $n=3$ and 
$$ G(\hat{\omega}) = G_3(\hat{\omega}) =
\frac{1}{2} \left(\frac{\hat{\omega}}{1-\hat{\omega}} +
\sqrt{
\left(\frac{\hat{\omega}}{1-\hat{\omega}} \right)^2+
 \frac{4 \hat{\omega}}{1-\hat{\omega}} }
 \right)
.
$$
was obtained in  \cite{mche}. Note that in an earlier paper \cite{R} a lot of geometric observation concerning the  behavior of the corresponding solutions of the system \eqref{s4} was done. In the present paper we want to carry out a deeper analysis of the geometry of the best approximations vectors in  the case  $n=3$ and introduce stronger new lower bounds for the ratio $\frac{\omega}{\hat{\omega}}$ under further condition.

\section{Patterns of best approximation vectors}

In the present paper it is convenient to deal with the best approximations with respect to the Euclidean norm. Let $|\pmb{x}|$ be the Euclidean norm of the vector $\pmb{x}$ in $\mathbb{R}^n$. 

Suppose that $\pmb{\alpha} \in \mathbb{R}^n\setminus \mathbb{Q}^n$.
We deal with the infinite  sequence of  best approximation vectors 
$$\pmb{z}_\nu = (q_\nu,\pmb{a}_\nu)= (q_\nu, a_{1,\nu},\ldots,a_{n,\nu})\in \mathbb{Z}^{n+1}, 
\,\,\,\,\,
\pmb{a}_\nu =(a_{1,\nu},\ldots,a_{n,\nu})\in \mathbb{Z}^{n},\,\,\,\,\,
\nu \ge 1$$

such that
$$q_1<q_2< \cdots <q_\nu<q_{\nu+1}< \cdots$$
form the increasing sequence of common denominators and the remainders of approximation
$$\xi_\nu =|q_\nu\alpha -\pmb{a}_\nu|$$ satisfy
$$\xi_1>\xi_2>\cdots>\xi_\nu>\xi_{\nu+1}> \cdots,$$
and
$$\xi_\nu = \min_{q\in \mathbb{Z}_+,\,\, q\le q_\nu; \,\, \pmb{a}\in \mathbb{Z}^n}
\,\,\,\,\,\,
|q\pmb{\alpha}-\pmb{a}|.$$

For approximations to one number, from continued fractions' theory we know that 
$$
\left|
\begin{array}{cc}
q_\nu &a_{1,\nu}
\cr
q_{\nu+1} &a_{1,\nu+1}
\end{array}\right|
= (-1)^\nu,\,\,\,\forall \, \nu \ge 1
$$
and in particular
each couple of consecutive  best approximations
$\pmb{z}_\nu, \pmb{z}_{\nu+1}$ consists of two linearly independent vectors.

Recall that the exponents $\omega(\pmb{\alpha})$ and $\hat{\omega}(\pmb{\alpha}) $ can be expressed in terms of best approximations by the formulas
\begin{equation}\label{laast}
\omega(\pmb{\alpha}) =-\liminf_{\nu \to \infty} \frac{\log \, \xi_\nu}{\log\, q_\nu},\,\,\,
\hat{\omega}(\pmb{\alpha}) = -\limsup_{\nu \to \infty} \frac{\log \, \xi_\nu}{\log\, q_{\nu+1}},\,\,\,
\frac{\omega(\pmb{\alpha})}{\hat{\omega}(\pmb{\alpha})} \ge  \limsup_{\nu \to\infty}\, 
\frac{\log\,q_{\nu+1}}{\log\, q_\nu}
\end{equation}
(for the last inequality here see our paper 
\cite{MaMo}, this inequality will be of importance to complete the proof of Theorem 1 at the very end of Section 5).

In the case $n = 2$  when $1,\alpha_1, \alpha_2$ are linearly independent over $\mathbb{Q}$,  Jarn\'{\i}k proved that  there  exist infinitely many $\nu$ such that 
three consecutive vectors $ \pmb{z}_{\nu-1}, \pmb{z}_\nu, \pmb{z}_{\nu+1}$ are linearly independent.

In the present paper we deal with the case $n=3$. The behavior of the best approximation vectors in this case is a little bit more complicated. From  Jarn\'{\i}k it follows again that  there  exist infinitely many $\nu$ such that three consecutive vectors $ \pmb{z}_{\nu-1}, \pmb{z}_\nu, \pmb{z}_{\nu+1}$ are linearly independent. However it may happen that there is no quadruples   $ \pmb{z}_{\nu-1}, \pmb{z}_\nu, \pmb{z}_{\nu+1}, \pmb{z}_{\nu+2}$ consisting of four consecutive independent vectors. It is the \emph{phenomenon of degeneracy of dimension}, see \cite{Che,Moor,MLMS,Mosin,Schlei}. Nevertheless, for any $\nu_0$ the sequence $ \pmb{z}_{\nu_0+j}, \, j \in \mathbb{Z}_+$ generates the whole space $\mathbb{R}^4$.
   
   Let us consider $\nu$ such that   three consecutive vectors $ \pmb{z}_{\nu-1}, \pmb{z}_\nu, \pmb{z}_{\nu+1}$ are linearly independent and define 
   $$
   j = \min \{ \nu'>\nu:\,\,
    \text{vectors} \,\, \pmb{z}_{\nu'-1}, \pmb{z}_{\nu'}, \pmb{z}_{\nu'+1}\,\,\text{are linearly independent}\}.
    $$ 
    Now we deal with the collection
     \begin{equation}\label{first}
    \pmb{z}_{\nu-1}, \pmb{z}_\nu, \pmb{z}_{\nu+1},
    \pmb{z}_{j-1}, \pmb{z}_j, \pmb{z}_{j+1},\,\,\,\, \nu< j
    \end{equation}
    of two successive triples 
    $ \pmb{z}_{\nu-1}, \pmb{z}_\nu, \pmb{z}_{\nu+1}$
    and
    $ \pmb{z}_{j-1}, \pmb{z}_j, \pmb{z}_{j+1}$ of independent vectors.
    We distinguish two cases.
    
    It may happen that  triples   $ \pmb{z}_{\nu-1}, \pmb{z}_\nu, \pmb{z}_{\nu+1}$ and $ \pmb{z}_{j-1}, \pmb{z}_j, \pmb{z}_{j+1}$ generate the same three-dimensional subspace of $\mathbb{R}^4$, that is
    $$
    \langle \pmb{z}_{\nu-1}, \pmb{z}_\nu, \pmb{z}_{\nu+1} \rangle_\mathbb{R} =  \langle \pmb{z}_{j-1}, \pmb{z}_j, \pmb{z}_{j+1}
    \rangle_\mathbb{R}.
    $$
    In this case we say the collection of vectors \eqref{first} has \emph{pattern} $A$.
    
    In the case when 
   triples   $ \pmb{z}_{\nu-1}, \pmb{z}_\nu, \pmb{z}_{\nu+1}$ and $ \pmb{z}_{j-1}, \pmb{z}_j, \pmb{z}_{j+1}$ generate different three-dimensional subspaces in $\mathbb{R}^4$, that is
    $$
    \langle \pmb{z}_{\nu-1}, \pmb{z}_\nu, \pmb{z}_{\nu+1} \rangle_\mathbb{R} \neq  \langle \pmb{z}_{j-1}, \pmb{z}_j, \pmb{z}_{j+1}
    \rangle_\mathbb{R}
    ,$$
    and so 
        $$
    \langle \pmb{z}_{\nu-1}, \pmb{z}_{j-1}, \pmb{z}_{j}, \pmb{z}_{j+1} \rangle_\mathbb{R} =
    \langle \pmb{z}_{\nu-1}, \pmb{z}_\nu, \pmb{z}_{\nu+1}, \pmb{z}_{j+1} \rangle_\mathbb{R}
    =\mathbb{R}^4,
    $$
    we say the collection of vectors  \eqref{first} has \emph{pattern} $B$.
     
    We should note that for each collection \eqref{first} there exists the next collection of the form \eqref{first}
      \begin{equation}\label{next}
    \pmb{z}_{j-1}, \pmb{z}_j, \pmb{z}_{j+1},
    \pmb{z}_{l-1}, \pmb{z}_l, \pmb{z}_{l+1},\,\,\,\, j<l
    \end{equation}
    where the last triple  in \eqref{first} coincide with the first triple in \eqref{next}. So, by encoding the collections \eqref{first} by letters $A$ and $B$ accordingly, for each vector $\pmb{\alpha}\in \mathbb{R}^3$ with 
    $1,\alpha_1,\alpha_2,\alpha_3$
    linearly independent over $\mathbb{Q}$  we define an infinite word
    \begin{equation}\label{word}
    \frak{W} (\pmb{\alpha}) =
    C_1 C_2 C_3 \ldots \, , \,\,\,\, C_j \in \{ A,B\}.
    \end{equation}
    From the discussion above we see that when
      $1,\alpha_1,\alpha_2,\alpha_3$ are 
    linearly independent over $\mathbb{Q}$,
    the word $
        \frak{W} (\pmb{\alpha}) $ is infinite and contains infinitely many letters $B$.

        In the present paper, we show how we can determine a sharper lower bound for the ratio of the exponents
        $\frac{\omega (\pmb{\alpha})}{\hat{\omega}(\pmb{\alpha})}$ in terms of the word 
        $
        \frak{W} (\pmb{\alpha}) $.
 
\section{Parameters}

 We  deal with parameters
\begin{equation}\label{pere}
\frac{1}{3}<\lambda < 1
\,\,\,\,
\text{and}
\,\,\,\, 0<
\theta= \theta(\lambda) = \frac{1-\lambda}{\lambda}<2.
\end{equation}
For $ k \ge 0$ we define
$$
  S_k = \sum_{j=0}^{k-1}\theta^j.
  $$
  So
  \begin{equation}\label{hh2}
    S_k
  \,=\,
  \begin{cases}
  k,\,\,\,\,\,\,\,\,\,\, \text{if}\,\,\,\,\,\theta = 1\cr
  \frac{\theta^k-1}{\theta-1},\,\,
  \text{if}\,\,\,\,\,
  \theta\neq 1
  \end{cases}
  \,\,\,\,\,\text{and}
  \,\,\,\,\,\,\,\,\,\,\,\,\,\,
  S_k
  \,\sim\,
  \begin{cases}
  \frac{1}{1-\theta},\,\,\,\,\, 
  \text{if}\,\,\,\,\,
  \theta<1 \cr
  k,\,\,\,\,\,\,\,\,\,\, \text{if}\,\,\,\,\,\theta = 1\cr
  \frac{\theta^k}{\theta-1},\,\,\,\,\,
  \text{if}\,\,\,\,\,
  \theta>1
  \end{cases}
  \,\,\,\,\,\,
  \text{when}\,\,\,\, 
  k \to \infty.
  \end{equation}

We consider a real variable $g$.  For $ k \ge1$ we define the quantities
$$
M_k = 
\theta S_k +\theta^{k+1} =\theta S_{k+1},
\,\,\,\,\,
N_k =2\theta^k+\frac{S_k}{\lambda}+\theta S_{k-1},
\,\,\,\,\,
P_k = 
\frac{S_{k-1}}{\lambda} +\theta^{k-1}
$$
and the polynomials
$$
R_k(g)= 
R_k(g,\lambda) = M_k g^2-N_k g +P_k.
$$
We are interested in the roots of polynomials $R_k(g)$.
For example,  
\begin{equation}\label{hhR1}
R_1(g,\lambda)  = (\theta+\theta^2)g^2 -\left(2\theta +\frac{1}{\lambda}\right) g +1,
\end{equation}
and the equation  $
R_1 (g,\lambda) = 0
$
may be written in the  form
\begin{equation}\label{xRoy}
R_1^* (g,\lambda)=
(1-\lambda) g^2 - (3-2\lambda)\lambda g +\lambda^2 = 0.
\end{equation}
It is clear that the equation (\ref{hhR1}) has a unique root $\frak{g}_1>1$ for which we have
\begin{equation}\label{hhqq}
\frak{g}_1>G(\lambda).
\end{equation}
Indeed,  the quantity
$G(\lambda)>1 $ satisfies the equation $(1-\lambda) g^2 -\lambda g -\lambda =0 $. 
As
$$
(1-\lambda) \left(
\frac{1+\lambda}{2(1-\lambda)}
\right)
^2 -\lambda \left(
\frac{1+\lambda}{2(1-\lambda)}
\right) -\lambda = \frac{1-3\lambda}{4} <0
,$$
we see that $ G(\lambda ) > 
\frac{1+\lambda}{2(1-\lambda)}>1
.$
Now
$$
R_1^*(G(\lambda ),\lambda) = \lambda G(\lambda ) +\lambda   - (3-2\lambda)\lambda G(\lambda )  +\lambda^2 = \lambda 
\left(
1+\lambda - 2(1-\lambda) G(\lambda)
\right) <0
$$
 and the existence of the root $\frak{g}_1$ satisfying  (\ref{hhqq}) follows.

 Here we should note  that $\theta = \theta(\lambda)$ depends on $\lambda$  and that the equation
$$
R_1(\theta(\lambda),\lambda) =0
$$
 has a root
$$\lambda_*=
\frac{1}{2}
 \left(2+\sqrt{5}-\sqrt{7+2\sqrt{5}}\right) = 0.42451^+
 $$
 which appeared in the problem of uniform simultaneous approximation to a real number, its square and its cube in a famous paper  \cite{R} by D. Roy.

\begin{lemma}\label{Lemma1}
 For any $ k\ge 1$ and for any $\lambda \in \left(\frac{1}{3}, 1\right)$  the equation
\begin{equation}\label{fu3}
R_k(g,\lambda)=0
\end{equation}
in the interval $\left( 1,\frac{2}{\theta}\right)$
has a unique root $\frak{g}_k = \frak{g}_k (\lambda)$.
This root has lower bound
\begin{equation}\label{GG}
 \frak{g}_k (\lambda) > G(\lambda).
\end{equation}
Moreover, for a fixed $\lambda$ one has 
\begin{equation}\label{fu4}
 \frak{g}_k (\lambda)< \frak{g}_{k+1} (\lambda)
 \,\,\,\,\,
 \text{and}
 \,\,\,\,\,
 \overline{\frak{g}}(\lambda)=
 \lim_{k\to \infty}  \frak{g}_k (\lambda) =
 \begin{cases}
 \frac{1}{1-\lambda},\,\,\,\text{if}\,\,\, \lambda \in \left[\frac{1}{2}, 1\right), \cr
 \frac{2}{\theta},\,\,\,\text{if}\,\,\, \lambda \in \left(\frac{1}{3},\frac{1}{2}\right].
 \end{cases}
\end{equation}
\end{lemma}
 
 \begin{proof} We  use the explicit form (\ref{hhR1}) of $R_1(g,\lambda)$  to see that
$$
R_1\left(\frac{1}{\theta}\right) = -1,\,\,\,\,\,
\text{and}
\,\,\,\,\,
R_1\left(\frac{2}{\theta}\right) = 
\frac{3\lambda -1}{1-\lambda} >0
$$
because of $ \lambda \in \left(\frac{1}{3},1\right)$. But for every $k\ge 1$ we have
\begin{equation}\label{pr}
R_{k+1}(g) -R_{k}(g)= \theta^k (\theta g -1)(\theta g -2) ,
\end{equation}
and so
$$
R_k\left(\frac{1}{\theta}\right) = R_1\left(\frac{1}{\theta}\right) <0,\,\,\,\,\,
\text{and}
\,\,\,\,\,
R_k\left(\frac{2}{\theta}\right) =  R_1\left(\frac{2}{\theta}\right)   >0.
$$
Moreover
$$
R_k(1)  = \theta^{k} (\theta- 2)< 0
$$
because of $\lambda >\frac{1}{3}$.
So for every $k\ge 1$, there exists a unique root $\frak{g}_k $ satisfying
 
\begin{equation}\label{prw}
\max\left\{1, \frac{1}{\theta}\right\} < \frak{g}_k < \frac{2}{\theta},
\end{equation}
as it was claimed in Lemma 1. Substituting $\frak{g}_k$ into (\ref{pr}) and taking into account (\ref{prw}) we see that 
$$
R_{k+1}(\frak{g}_k)<0,
$$
and this shows that 
$\frak{g}_{k+1}>\frak{g}_k$. 
Now (\ref{GG}) follows from (\ref{hhqq}) and monotonicity. 
The limit condition from (\ref{fu4})  follows from (\ref{hh2}), because the limit value $\overline{\frak{g}}$ satisfies equation
$$
(1-\lambda)g^2 -(2-\lambda)g + 1=
((1-\lambda)g-1)(g-1) = 0
\,\,\,\,\,\text{for}\,\,\,\,\, \lambda \in \left[\frac{1}{2}, 1\right),
$$
or
$$
\theta^2g^2 - 3\theta g + 2 =  (\theta g-1)(\theta g - 2) =0
\,\,\,\,\,\text{for}\,\,\,\,\, \lambda \in \left(\frac{1}{3}, \frac{1}{2}\right].
$$
Lemma 1 is proven.
\end{proof}
 
\begin{remark} 
From the inequalities $ G(\lambda )>\frac{\lambda}{1-\lambda}$ and \eqref{GG},\eqref{fu4}  we see that for any $k\ge 1$
\begin{equation}\label{hh1}
\frac{1}{\theta}=
 \frac{\lambda}{1-\lambda}< \frak{g}_k (\lambda)< \frac{1}{1-\lambda}
 .
\end{equation}
\end{remark}

\section{ Results}

We consider the quantity 
$$
k = k(\pmb{\alpha}) =  \sup \left\{ j:\,\, \text{subword}\,\, \underbrace{A\ldots A}_{j\, \text{times}}B\,\,
\text{occurs infinitely many times in \eqref{word}}\right\}
$$
which may be finite or infinite.

First of all, let us discuss the case $ k(\pmb{\alpha}) =0$. In this case the word  $\frak{W} (\pmb{\alpha})$ has only finitely many letters $A$ and we do not improve on the bound $\frac{\omega(\alpha)}{\hat{\omega}(\alpha)} \ge G(\hat{\omega})$
which is known to be optimal in this case (see \cite{R1} and the explanation from \cite{SS_MJCNT}). 
The extremal cases $\hat{\omega}(\pmb{\alpha}) =\frac{1}{3}$ and $\hat{\omega}(\pmb{\alpha}) =1$, where we have respectively
$\omega(\pmb{\alpha}) \ge \frac{1}{3}$ or  $\omega(\pmb{\alpha}) =\infty$ are eluded from our further considerations.

\begin{maintheorem}\label{Thm1}
Suppose that  numbers $ 1,\alpha_1,\alpha_2,\alpha_3$ are linearly independent over $\mathbb{Q}$ and suppose that 
  $k(\pmb{\alpha}) = k\ge 1$ is finite. Suppose that $ \frac{1}{3} <\hat{\omega}(\alpha)<1$.
  Then
$$
\frac{\omega(\pmb{\alpha})}{\hat{\omega}(\pmb{\alpha})} \ge \frak{g}_k(\hat{\omega}(\alpha)),
$$
where $\frak{g}_k (\cdot)\ge 1$ is the root of  \eqref{fu3} defined in Lemma \ref{Lemma1}.
In the case   $k(\pmb{\alpha}) =\infty$  we have the inequality
$$
\frac{\omega(\pmb{\alpha})}{\hat{\omega}(\pmb{\alpha})} \ge \overline{\frak{g}}(\hat{\omega}(\alpha)).
$$

\end{maintheorem}

In the next statement we show the optimality of the bound from Theorem \ref{Thm1}.

\begin{maintheorem}\label{Thm2}
For any $\lambda \in \left(\frac{1}{3}, 1\right) $
there exists $\pmb{\alpha}\in \mathbb{R}^3$  such that the sequence of patterns 
\eqref{word} for 
 $\pmb{\alpha}$ is an ultimately
 periodic sequence  with the period
 $$
 \underbrace{AA\ldots A}_{k \, \text{times}} B
 $$
and
\begin{equation}\label{T2}
\hat{\omega}(\pmb{\alpha}) = \lambda,
\,\,\,\,\,
\,\,\,\,\,
\frac{\omega(\pmb{\alpha})}{\hat{\omega}(\pmb{\alpha})} = \frak{g}_k (\lambda).
\end{equation}
 \end{maintheorem}
 
 \begin{remark} 
 The vector $\pmb{\alpha}$ constructed in Theorem 2 has
 $1,\alpha_1,\alpha_2,\alpha_3$ linearly independent over $\mathbb{Q}$, because pattern $B$ occurs in the word \eqref{word} infinitely often.
 \end{remark}
 
\begin{remark} 
It is possible  to show that any sequence from two  letters $A,B$ with infinitely many $B$'s can occur as the word
$\frak{W}(\pmb{\alpha})$ for a certain $\pmb{\alpha}$ with linearly independent components.
Moreover, for any $\lambda  \in \left(\frac{1}{3}, 1\right)  $ and  for any word $\frak{W}$  with infinitely many $B$'s, it is possible to construct
$\pmb{\alpha}$ with linearly independent components such that 
$\frak{W}(\pmb{\alpha})=\frak{W}$  and for 
 $ k =k(\pmb{\alpha})$ (determined by the word $\frak{W}$)
the equalities \eqref{T2} are valid.
\end{remark}

We  give a proof of Theorem \ref{Thm1} in Section \ref{5} below. As for Theorem \ref{Thm2}, it is proved by a certain inductive construction of the limit vector $\pmb{\alpha}$. The idea of such a construction is quite easy, one should construct a sequence of integer vectors which 
form two- and thee-dimensional subspaces with desired properties and which turn out to form the sequence of all the best approximation vectors to $\pmb{\alpha}$.
As usual such constructions are rather cumbersome (see, for example 
\cite{MX,MY}).
Here  we give a detailed proof of the case $k=1$ in Sections 6 - 8.  The proof for the case $k>1$ uses similar argument. The difference is that  in the inductive process described below we need to repeat $k$ times {\bf Stage 1}. All necessary comments for the case $k>1$ are given in Section \ref{9}.\\

\begin{remark}
Note that it should be possible and interesting to use the theory of parametric geometry of numbers to prove Theorem \ref{Thm1}. However, it seems to the authors that it cannot be used directly to prove Theorem  \ref{Thm2}. Indeed, given a system, it does not seem possible to associate it to a parameter $k$. It may be at bounded distance to successive minima functions of points with different parameter $k$.
\end{remark}

\section{Proof of Theorem 1}\label{5}

 Let us consider $\lambda$ satisfying 
  $\frac{1}{3}< \lambda <\hat{\omega} (\pmb{\alpha}) $. Then for all $\nu$ large enough one has
\begin{equation}\label{sing}
\xi_\nu \le q_{\nu+1}^{-\lambda}.
\end{equation}

  Let
 \begin{equation}\label{triples}
     \pmb{z}_{\nu_k-1}, \pmb{z}_{\nu_k}, \pmb{z}_{\nu_k+1};\,\,\,
     \pmb{z}_{\nu_{k-1}-1}, \pmb{z}_{\nu_{k-1}}, \pmb{z}_{\nu_{k-1}+1};\,\,\, \ldots \,\,\, ; \,\,\,
     \pmb{z}_{\nu_1-1}, \pmb{z}_{\nu_1}, \pmb{z}_{\nu_1+1};\,\,\,
      \pmb{z}_{\nu_0-1}, \pmb{z}_{\nu_0}, \pmb{z}_{\nu_0+1};\,\,\,
     \pmb{z}_{l-1}, \pmb{z}_{l}, \pmb{z}_{l+1}
 \end{equation}
 be the triples of successive best approximation vectors corresponding to the subword 
 $ \underbrace{A\ldots A}_{k\, \text{times}}B$,
 so
 $$
 \nu_k<\nu_{k-1}<\cdots<\nu_1<\nu_0< l
 $$
 and each successive  triples
  $$
     \pmb{z}_{\nu_i-1}, \pmb{z}_{\nu_i}, \pmb{z}_{\nu_i+1};\,\,\,
     \pmb{z}_{\nu_{i-1}-1}, \pmb{z}_{\nu_{i-1}}, \pmb{z}_{\nu_{i-1}+1},\,\,\,\,\, 1\le i \le k
     $$
     determine pattern $A$, while the successive  triples
     $$
        \pmb{z}_{\nu_0-1}, \pmb{z}_{\nu_0}, \pmb{z}_{\nu_0+1};\,\,\,
     \pmb{z}_{l-1}, \pmb{z}_{l}, \pmb{z}_{l+1}
 $$
 determine pattern $B$.

 Consider the two-dimensional lattice
 $$
   \Lambda = \langle \pmb{z}_{\nu_0},\pmb{z}_{\nu_{0}+1}\rangle_{\mathbb{Z}}  =
   \langle \pmb{z}_{{l}-1},\pmb{z}_{{l}}\rangle_{\mathbb{Z}}
   $$
   and the three-dimensional lattices
    $$
    \Gamma_{i} = \langle \pmb{z}_{\nu_i-1}, \pmb{z}_{\nu_i},\pmb{z}_{\nu_{i}+1}\rangle_{\mathbb{Z}} ,\,\, i = 0,\ldots,k
      \,\,\,\,\,  \text{and}\,\,\,\, \,\, 
  \Gamma = 
   \langle \pmb{z}_{{l}-1},\pmb{z}_{{l}}\, \pmb{z}_{l+1}\rangle_{\mathbb{Z}}.
   $$
    Consider the  two-dimensional subspace
  $  \frak{L}= {\rm span}\, \Lambda $.
    As for three-dimensional subspaces,
   from the definitions of patterns $A$ and $B$ we see that 
   $$
   \frak{G}^*=
   {\rm span} \,\Gamma_0 =
      {\rm span} \,\Gamma_1= \cdots
      =
         {\rm span} \,\Gamma_k
   $$
   is the same three-dimensional subspace,
   $$
      \frak{G} =    {\rm span} \,\Gamma,\,\,\,\,
      {\rm span} \,(\Gamma_{i} \cup \Gamma) = \mathbb{R}^4,\,\,\,\,
      0\le i \le k,
      $$       
      and
 $$
    ({\rm span} \,\Gamma_{i} )     \cap 
      ( {\rm span} \,\Gamma)  
      = \frak{G}^* \cap \frak{G} = \frak{L},\,\,\,\, 0\le i \le k.
      $$
      Moreover, vectors  
      \begin{equation}\label{ind}
      \pmb{z}_{\nu_0-1}
      ,
        \pmb{z}_{l-1},
        ,
          \pmb{z}_{l}
          ,
            \pmb{z}_{l+1}
            \end{equation}
            are independent.
      As  vectors \eqref{ind} are independent vectors in $\mathbb{R}^4$, their coordinates form a non-zero integer determinant and
      $$
         1\ll 
        \xi_{\nu_0-1}\xi_{l-1}\xi_{l} q_{l+1},
        $$
        or taking into account \eqref{sing} and the definition  of $\theta$ from \eqref{pere},
        \begin{equation}\label{tw2}
        q_{\nu_0}\ll q_{l}^{-1} q_{l+1}^{\theta}
 .
        \end{equation}
            By standard argument we have the bounds
      \begin{equation}\label{tw}
      {\rm det}\, \Lambda \asymp \xi_{l-1}q_l
      \end{equation}
      and
         \begin{equation}\label{tw1}
      {\rm det}\, \Gamma_{i} \ll \xi_{\nu_i-1}\xi_{\nu_i}q_{\nu_i+1},\,\, 0\le i \le k,\,\,\,\,\,
      {\rm det}\, \Gamma \ll \xi_{l-1}\xi_{l}q_{l+1}.
      \end{equation}
      A very detailed proof of (\ref{tw}) and (\ref{tw1}) is given for example in  Lemma 2 from \cite{MaMo}. We should recall that the signs $\asymp$ and $\ll$ in 
      (\ref{tw}) and (\ref{tw1}) contain {\it absolute} constants.
      From (\ref{sing}) we  deduce the inequality
      \begin{equation}\label{44444}
      q_{j+1}\gg q_j^{\frac{1}{\theta}},\,\,\,\,\, j = l, \nu_0,\nu_1,...,\nu_k
      \end{equation}
      which is trivial for $ \lambda < \frac{1}{2}$.
        Now we consider the rational subspaces $\frak{G}, \frak{G}^*$ and $\frak{L}$ defined above.  For their heights
           \begin{equation}\label{tw4}
           H (\frak{G}^*) \le    {\rm det}\, \Gamma_k,\,\,\,\,
             H (\frak{G}) \le    {\rm det}\, \Gamma,
             \,\,\,\,
             H(\frak{L}) \asymp   {\rm det}\, \Lambda
             \end{equation}
             by Schmidt's inequality 
             (see \cite{SchHeights}, Lemma 8A)
             we have the  bound
                   \begin{equation}\label{tw5}
             1\ll   H (\frak{G}^*) \cdot \frac{ H (\frak{G})  }{  H(\frak{L})}
           \end{equation}
Now \eqref{tw4},\eqref{tw5} together with \eqref{tw},\eqref{tw1}  gives us
 $$
 1\ll   \xi_{\nu_k-1}\xi_{\nu_k}q_{\nu_k+1} \cdot \frac{\xi_{l}q_{l+1}}{q_{l}},     
  $$
  or by \eqref{sing},
            \begin{equation}\label{tw6}
            1\ll q_{\nu_k}^{-\lambda} q_{\nu_k+1}^{1-\lambda} \cdot q_l^{-1} q_{l+1}^{1-\lambda},
            \end{equation} 

\begin{lemma}\label{Lemma2}
  Suppose that $g$ satisfies 
 \begin{equation}\label{55555}
\frac{1}{\theta}=
 \frac{\lambda}{1-\lambda}< g< \frac{1}{1-\lambda}
\end{equation}
and define
 \begin{equation}\label{66666}
\sigma = \sigma (g,\lambda) =\frac{\frac{1}{1-\lambda}-g}{ g - \frac{1}{\theta}} >0
\end{equation}
For $ k\ge 1$ we consider subword   $ \underbrace{A\ldots A}_{k\, \text{times}}B$ and use the notation defined before.

 Suppose that  
 \begin{equation}\label{gg}
 q_{l+1}\le q_l^g.
 \end{equation}
  Then we have 
  \begin{equation}\label{biblio}
   \theta^k >\sigma S_k
   \end{equation}
   and
 \begin{equation}\label{fu2}
 q_{\nu_k+1}\gg  q_{\nu_k}^{f_k},
 \,\,\,\,\,\text{
 where }
  \,\,\,\,\,
f_k  = \frac{ \theta^{k-1}- \sigma S_{k-1} }{\theta^{k}- \sigma S_{k} }
.
\end{equation}

   \end{lemma}
 
 \begin{proof}
 We proceed by induction in $k$. 
 
 First we verify the base of induction for $k=1$.
 
  From $\nu_1<\nu_0$ and \eqref{tw2},\eqref{gg} we deduce
 \begin{equation}\label{AA}
 q_{\nu_1+1}\le q_{\nu_0}\le q_l^{-1+\theta g}.
 \end{equation}
Using (\ref{tw6}) and (\ref{gg}), one gets
$$
1\ll q_{\nu_1}^{-\lambda} q_{\nu_1+1}^{1-\lambda} q_l^{-1+(1-\lambda)g}.
$$
As the exponent of $q_l$ in this expression is negative, one may use (\ref{AA}) to eliminate $q_l$ from this estimate. Since 
$$
1-(1-\lambda) g  =\sigma\lambda (\theta g - 1)
$$
by (\ref{66666}), this yields 
$$
1\ll q_{\nu_1}^{-\lambda}q_{\nu_1+1}^{1-\lambda -\sigma \lambda}
$$
and so
$$
q_{\nu_1+1}^{\theta-\sigma} \gg q_{\nu_1}.
$$
As $l$ and thus $\nu_1$ can be chosen  arbitrary large, this implies that $\theta> \sigma$ and
$q_{\nu_1+1}\gg q_{\nu_1}^{f_1}$.
So we have checked (\ref{biblio}) and (\ref{fu2}) for $k=1$.

 \vskip+0.3cm
 
 The induction step can be treated in the same way. Using 
  (\ref{tw6}), (\ref{gg}) and (\ref{AA}), one gets
  $$
  1\ll q_{\nu_{k+1}}^{-\lambda} q_{\nu_{k+1}+1}^{1-\lambda}q_{\nu_0}^{-\sigma\lambda}.
  $$
   By inductive assumption
 $$
 q_{\nu_i+1}\ge q_{ \nu_i}^{f_i},\,\,\,\,\,
 1\le i \le k $$
 and so as
 $$
 q_{\nu_i+1}\le q_{\nu_{i-1}},\,\,\,1\le i \le k
 ,
 $$
 we get
 \begin{equation}\label{CC}
 q_{\nu_k}\le q_{\nu_0 }^{\frac{1}{\prod_{j=1}^k f_j}}.
 \end{equation}
  Also by induction we have $\theta^k > \sigma S_k$. So $f_1, \ldots ,f_k$ are defined through (\ref{fu2})
  and the estimate (\ref{CC}) simplifies to
  $$
  q_{\nu_0} \gg q_{\nu_k}^{\frac{1}{\theta^k - \sigma S_k}}.
  $$
  We take into account the inequality $ q_{\nu_{k+1}+1}\le q_{\nu_k}$ and deduce that
 $$
 q_{\nu_{k+1}}\ll q_{\nu_{k+1}+1}^{\theta-\frac{\sigma}{\theta^k-\sigma S_{k+1}}}
 = q_{\nu_{k+1}+1}^{\frac{\theta^{k+1}-\sigma S_{k+1}}{\theta^k -\sigma S^k}}.
 $$
 As $l$ and thus $\nu_{k+1}$ can be chosen arbitrary large, this implies that $\theta^{k+1} > \sigma S_{k+1}$ and $q_{\nu_{k+1}+1}\gg q_{\nu_{k+1}}^{f_{k+1}}$.
 
  \end{proof}

 \begin{remark}
 The functions $f_k$ defined in (\ref{fu2})  satisfy the nice  recursive formulas
\begin{equation}\label{india}
f_k = \frac{1}{\frac{1}{\lambda} - \theta f_{k-1}}, \,\,\,\,\, k \ge 2,
\end{equation}
so
$$
f_k =
\frac{1}{\displaystyle{\frac{1}{\lambda}-\frac{\theta}{\displaystyle{\frac{1}{\lambda}-
\frac{\theta}{\displaystyle{\frac{1}{\lambda} -\dots -
\frac{\theta}{
\displaystyle{ \frac{1}{\lambda} -\frac{\theta}{\theta-\sigma}}} }}}}}} 
$$
(here $\frac{1}{\lambda}$ occurs $k-1$ times).
 \end{remark}

 Now we finish with the proof of Theorem \ref{Thm1}. 
 We can suppose that \eqref{gg} holds with $g <\frac{1}{1-\lambda}$,
  otherwise $ \frac{\omega}{\hat{\omega}} \ge \frac{1}{1-\lambda}$ and there is nothing to prove.
By Lemma \ref{Lemma2} we  deduce \eqref{fu2}.

Now we use  (\ref{pere}), (\ref{66666}) and (\ref{fu2}) to see that  equation (\ref{fu3}) which determines the values of $\frak{g}_k$
can be written as 
$$
f_k(g,\lambda)=g.
$$
 This means  that
 $$
 \text{either}\,\,\,\,
 q_{l+1}\gg q_l^{\frak{g}_k}\,\,\,\,\,
 \text{or}
 \,\,\,\,\,
  q_{\nu_k+1}\gg q_{\nu_k}^{\frak{g}_k}.
  $$
  and \eqref{laast} gives the bound $ \frac{\omega}{\hat{\omega}} \ge \frak{g}_k$ and Theorem \ref{Thm1} is proved. \qed

\section{Inductive construction for Theorem 2: $k=1$}

 To deal with the proof of Theorem 2 in the case $k=1$, we need more parameters.
 Define
\begin{equation}\label{g11} 
\frak{g}_{1,0} = \frac{1}{\theta \frak{g}_1-1 }
 = \frac{\frak{g}_1-\lambda}{(1-\lambda) \frak{g}_1 }
=
\frac{1-\lambda-\frac{\lambda}{\frak{g}_1}}
{1-(1-\lambda)\frak{g}_1}
\end{equation}
where equalities follow from the equation \eqref{xRoy} which determines the root $\frak{g}_1$.
 
Note that for $\lambda \in \left(\frac{1}{3},1\right)$ one has
\begin{equation}\label{g1}
1\le 
\frak{g}_{1,0} \le \frak{g}_1.
\end{equation}
Indeed, from Lemma \ref{Lemma1} we see that $ \frak{g}_1 < \frac{2}{\theta} $,
and this together with the definition (\ref{g11}) gives $\frak{g}_{1,0} >1$. To explain the second inequality from (\ref{g1}) we recall that  by Lemma \ref{Lemma1} one has $ \frak{g}_{1}\ge G(\lambda)$ where $G(\lambda)$ is the unique positive root of the equation $ \theta x^2 = x+1$. This means that 
  $ \theta \frak{g}_{1}^2 \ge \frak{g}_{1}+1$,
  or
    $ (\theta \frak{g}_{1} -1) \frak{g}_{1}\ge 1$,
   and the right  inequality from (\ref{g1}) follows.

We construct 
$\pmb{\alpha}\in  [-1,1]^3 \subset \mathbb{R}^3$
 such that  for any  $j$, its best approximation vectors $\pmb{z}_j = (q_j, a_{1,j},a_{2,j},a_{3,j})$
satisfy
 \begin{equation}\label{iv}
 \xi_{j} \asymp q_{j+1}^{-\lambda}
 \end{equation}
 and moreover for any $\nu$, the five successive best approximation vectors
\begin{equation}\label{collek}
\pmb{z}_{2\nu-1}, \pmb{z}_{2\nu}, \pmb{z}_{2\nu+1},\pmb{z}_{2\nu+2},\pmb{z}_{2\nu+3}
\end{equation}
satisfy the conditions ({\bf i}) - ({\bf vii}) below.
 For $\pmb{z}_j$ we consider a shortened vector $\underline{\pmb{z}}_j = (a_{1,j},a_{2,j}, a_{3,j})$.
We denote by
 $$
 \pmb{Z}_j = q_j^{-1} \underline{\pmb{z}}_j= \left(
 \frac{a_{1,j}}{q_j}, \frac{a_{2,j}}{q_j},  \frac{a_{3,j}}{q_j}\right)
 \in \mathbb{Q}^3
 $$
 the corresponding rational vector.
 
 As we are constructing  vector $\pmb{\alpha}$   from the cube
$ [0,1]^3$ we may assume that and the approximating vectors $\pmb{z}=(q,a_1,a_2, a_3)\in \mathbb{Z}^4$ from the proof below, including  the best approximation vectors,
satisfy the inequality
$$
\max_{j=1,2,3}\left|\alpha_j - \frac{a_j}{q}\right| \le 1
$$
and so
\begin{equation}\label{fri}
q\le |\pmb{z}| \le  4q.
\end{equation}

 We say that  vectors 
\begin{equation}\label{bas}
\pmb{u},\pmb{v},\pmb{w}\in \mathbb{Z}^4
\end{equation}
form a \emph{primitive triple} if
$$
\langle
\pmb{u},\pmb{v},\pmb{w}
\rangle_\mathbb{Z} =
\mathbb{Z}^4 \cap
\langle
\pmb{u},\pmb{v},\pmb{w}
\rangle_\mathbb{R} 
,
$$
that is  vectors \eqref{bas} form a basis of the three-dimensional lattice 
$\mathbb{Z}^4\cap 
\langle
\pmb{u},\pmb{v},\pmb{w}
\rangle_\mathbb{R} 
.$

 Recall that for a vector $\pmb{\zeta}$ in $\mathbb{R}^3$ or
 $\mathbb{R}^4$ we use the notation $ |\pmb{\zeta}|$ for the Euclidean norm. For two successive vectors $\pmb{z}_j$ and
 $\pmb{z}_{j+1}$ we consider the value
 $$
 \zeta_j 
= \left| \pmb{z}_j- \frac{q_j}{q_{j+1}}\pmb{z}_{j+1}\right|=
 q_j |\pmb{Z}_j - \pmb{Z}_{j+1}|.
 $$
  For three points $\pmb{Z}, \pmb{Z}'\neq \pmb{Z},\pmb{Z}''\neq \pmb{Z}\in \mathbb{R}^3$
 we define $ \angle  \pmb{Z}'\pmb{Z}\pmb{Z}''$ to be the (non-oriented) angle between non-zero vectors 
 $\pmb{Z}'-\pmb{Z}$ and  $\pmb{Z}''-\pmb{Z}$.

 \vskip+0.3cm
 
  The desired conditions for the best approximations vectors $\pmb{z}_j$ are as follows:\
 
  \vskip+0.3cm
 
\noindent
({\bf i})  vectors
$
\pmb{z}_{2\nu-1}, \pmb{z}_{2\nu}, \pmb{z}_{2\nu+1} $, as well as vectors
$
\pmb{z}_{2\nu+1}, \pmb{z}_{2\nu+2}, \pmb{z}_{2\nu+3} $
  form  primitive triples;

\noindent
({\bf ii}) with some integers $ a_\nu$ and $b_\nu$  we have
$\pmb{z}_{2\nu+2} = \pmb{z}_{2\nu-1} + a_\nu\pmb{z}_{2\nu}+ b_\nu\pmb{z}_{2\nu+1}$ and so 
$
  \pmb{z}_{2\nu}, \pmb{z}_{2\nu+1},\pmb{z}_{2\nu+2} $  form a primitive triple;

\noindent
({\bf iii})   vectors
$
\pmb{z}_{2\nu-1}, \pmb{z}_{2\nu}, \pmb{z}_{2\nu+1} , \pmb{z}_{2\nu+3}$  
(and so  vectors
$
\pmb{z}_{2\nu}, \pmb{z}_{2\nu+1}, \pmb{z}_{2\nu+2} , \pmb{z}_{2\nu+3}$) 
form a  basis of the lattice $\mathbb{Z}^4$;

\noindent
({\bf iv})  \,\, \,
 $q_{2j+1}\asymp q_{2j}^{\frak{g}_1},\,\,\,\,\,\,\,     \zeta_{2j} \asymp  q_{2j+1}^{-\lambda} \asymp q_{2j}^{-\lambda\frak{g}_1},\,\,\,\,\,\,\,\,\,\,\,\,\, j =\nu, \nu+1$;

\noindent
({\bf v})  \,\, \,
 $q_{2j+2}\asymp q_{2j+1}^{\frak{g}_{1,0}},\,\,\,
 \zeta_{2j+1} \asymp q_{2j+2}^{-\lambda} \asymp q_{2j+1}^{-\lambda\frak{g}_{1,0}}
 ,
 \,\,\,\, j =\nu-1, \nu,\,\,\,\,\text{and}\,\,\,\, \zeta_{2j+3} \asymp q_{2j+4}^{-\lambda} \asymp q_{2j+3}^{-\lambda\frak{g}_{1,0}}
 $;

\noindent
({\bf vi})  \,\,\,  the angle  between rational three-dimensional vectors
$\pmb{Z}_{2\nu+1} - \pmb{Z}_{2\nu+3}$ and $\pmb{Z}_{2\nu+2}- \pmb{Z}_{2\nu+3}$ is 
greater that $\frac{\pi}{4}$ and less than $\frac{3\pi}{4}$, that is
$$
\frac{\pi}{4}<
\angle \pmb{Z}_{2\nu+1}  \pmb{Z}_{2\nu+3}\pmb{Z}_{2\nu+2}<
\frac{3\pi}{4}
;
$$

\noindent
({\bf vii})  for the fundamental volume $\Delta$  of three-dimensional lattice
$$
\Gamma=  \langle\pmb{z}_{2\nu+1},\pmb{z}_{2\nu+2},\pmb{z}_{2\nu+3}
\rangle_\mathbb{Z}  = \mathbb{Z}^4 \cap  \langle\pmb{z}_{2\nu+1},\pmb{z}_{2\nu+2},\pmb{z}_{2\nu+3}
\rangle_\mathbb{R}
$$
we have 
\begin{equation}\label{wowo}
\Delta \asymp
 \zeta_{2\nu +1}\zeta_{2\nu+2} q_{2\nu +3}
 .
 \end{equation}
 
  \vskip+0.3cm 

\begin{remark}
Here in conditions ({\bf iv}),  ({\bf vi}), ({\bf vii}), in inequality (\ref{iv}) and in the sequel constants in signs
$\ll, \asymp$ may depend on $\lambda$ (and later on $k$) but they do not depend on $\nu$, that is for all values of $\nu$ the constants are the same. So the constants here do not affect the asymptotic behaviour of the ratios
$\frac{\log q_{j+1}}{\log q_j}$ and hence the exponents $\omega$ and $\hat{\omega}$.
\end{remark}

\vskip+0.3cm 
If we construct  $\pmb{\alpha}$ for which the sequence of the best approximation vectors $\pmb{z}_j$ satisfy the conditions above, we prove Theorem \ref{Thm2} for $ k=1$. Indeed, 
as
$$
\pmb{Z}_j \to \pmb{\alpha},\,\,\,\,\,\frac{\xi_j}{\zeta_j} \to 1\,\,\,\,\, \text{when}\,\,\,\,\, j \to \infty,
$$
from the conditions  ({\bf iv},{\bf v}) 
we see that $ q_{2j+1}\asymp q_{2j}^{\frak{g}_1}, q_{2j+2}\asymp q_{2j+1}^{\frak{g}_{1,0}}$  and 
\eqref{iv} hold for every $j$, and the inequality $ \frak{g}_{1,0}\le \frak{g}_1$ gives us
$\frac{\omega(\alpha)}{\hat{\omega}(\alpha)} = \frak{g}_1$.
Meanwhile, the conditions ({\bf i},{\bf ii}) and ({\bf iii})   show that the sequence of patterns for $\pmb{\alpha}$ is just the periodic sequence $ABABABABAB\ldots$ .

Vector $\pmb{\alpha}$ is constructed as a result of a certain inductive procedure and conditions ({\bf vi}, {\bf vii}) are necessary to proceed with inductive step.

%We need some further objects.
 %Suppose that  $\pmb{Z}_{j-2}, \pmb{Z}_{j-1}$ and $\pmb{Z}_j$ do not belong to the same {\color{blue} affine line in $\mathbb{R}^3$.
  % Then the points  $\pmb{Z}_{j-2}, \pmb{Z}_{j-1},\pmb{Z}_j$ determine a unique affine plane in $\mathbb{R}^3$.}
  Now we  define the \emph{angular neighbourhood}
  $$
 \mathcal{U}_j =
 \{\pmb{X}\in \mathbb{R}^3:\,\,\,\,
 |\pmb{X}-\pmb{Z}_j|\le  q_j^{-u_j}\,\,\,\,
 \text{and}\,\,\,\,
 \frac{\pi}{4}< \angle  \pmb{Z}_{j-1}\pmb{Z}_j\pmb{X}<
  \frac{\pi}{2}\},
 $$
 where
\begin{equation}\label{uji}
 u_j =
 \begin{cases}
 1+ \lambda \frak{g}_{1,0},\,\,\,\,\,\,\, j\,\, \text{odd},
 \cr
 1+ \lambda \frak{g}_1,\,\,\,\,\,\,\,\,\,\, j\,\, \text{even}.
 \end{cases}
\end{equation}
   For a point $\pmb{Y}\in \mathbb{R}^3$ we consider the ball
 $$\mathcal{W}_j (\pmb{Y}) =  \{\pmb{X}\in \mathbb{R}^3:\,\,\,\,
 |\pmb{X}-\pmb{Y}|\le \frac{1}{10} \cdot q_j^{-u_j}\}.$$
The neighborhoods under consideration are shown on the Fig. 1 below.
   
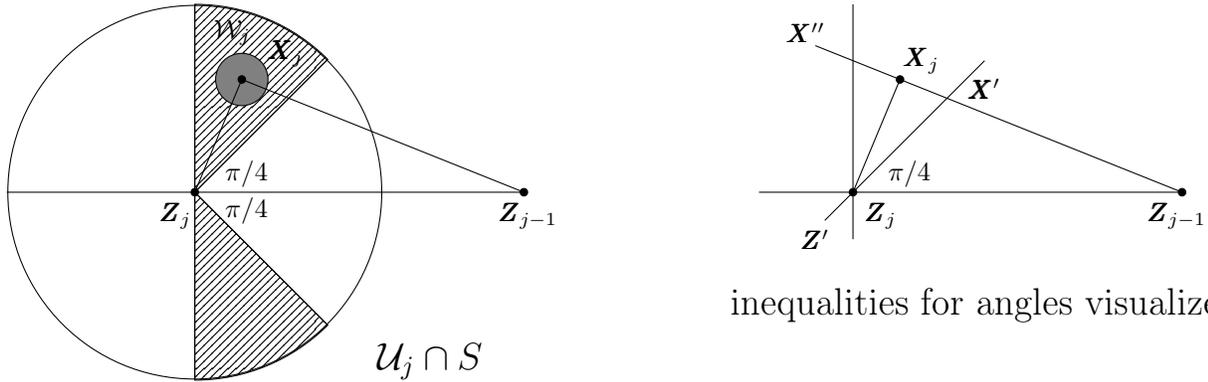
\begin{figure}[h]
  \centering

\begin{tikzpicture}[scale=1.25]

  \node[draw=black,fill=white,circle,inner sep=50pt] at (-5, 0) {};
  \node[draw=black,fill=black,circle,inner sep=1pt] at (-5, 0) {};
   
  \draw [pattern=north east lines, fill opacity=0.25] (-5,0) -- (-5,-2) arc (-90:-45:2) -- cycle;
  \draw [pattern=north east lines, fill opacity=0.25] (-5,0) -- (-5,2) arc (90:45:2) -- cycle;
  
  \node[draw=black,fill=black,circle,inner sep=1pt] at (-1.5, 0) {};
    
  \draw (-7,0) -- (-1.5,0);
  \draw (-5,-2) -- (-5,2);
  \draw (-5,0) -- (-3.6,1.4);
  \draw (-5,0) -- (-3.6,-1.4);
          
  \draw (-5.2, -0.25)   node { $\pmb{Z}_j$};  
  \draw (-1.45, -0.25)   node { $\pmb{Z}_{j-1}$};      
          
  \node[draw=black,fill=gray,circle,inner sep=7pt] at (-4.5,1.2) {};
           
  \draw (-4.5,1.2) -- (-1.5,0);
  \draw (-4.5,1.2) -- (-5,0);
                
  \node[draw=black,fill=black,circle,inner sep=1pt] at (-4.5,1.2) {};
                  
  \draw (-4.6, 1.7)   node { $\mathcal{W}_j$};      
  \draw (-4.03, 1.48)   node { $\pmb{X}_j$};   
  \draw (-4.45, 0.2)   node { $\pi/4$};   
  \draw (-4.45, -0.2)   node { $\pi/4$};  
  \draw (-2.5, -1.8)   node { \begin{Large}$\mathcal{U}_j\cap S$\end{Large}};

%%%%%%%%%%%%%%%%%%%%%%%%%%%%%%%%%%%%%%%

   \draw (2.73, 1.4)   node { $\pmb{X}_j$};  
   \draw (2.6, 0.2)   node { $\pi/4$}; 
   \draw (2.3, -0.25)   node { $\pmb{Z}_j$};  
   \draw (5.45, -0.25)   node { $\pmb{Z}_{j-1}$};   
   \node[draw=black,fill=black,circle,inner sep=1pt] at (2, 0) {}; 
   \node[draw=black,fill=black,circle,inner sep=1pt] at (5.5, 0) {}; 
   \draw (1,0) -- (5.5,0);
            
   \draw (1.6,1.56) -- (5.5,0);
   \draw (2.5,1.2) -- (2,0);
   \node[draw=black,fill=black,circle,inner sep=1pt] at (2.5,1.2) {};
   \draw (2,-0.5) -- (2,2);
   \draw (1.7,-0.3) -- (3.4,1.4);
   \draw (1.6, -0.5)   node { $\pmb{Z}'$};  
   \draw (1.5, 1.7)   node { $\pmb{X}''$};  
   \draw (3.4, 1.1)   node { $\pmb{X}'$}; 
    \draw (3.4, -1.2)   node { \begin{Large}\text{inequalities for angles visualized}\end{Large}}; 

\end{tikzpicture}

      \caption{Angular domain $\mathcal{U}_j$ and the neighborhood $\mathcal{W}_j$ inside, section by subspace $S$; inequality for angles
      $ \frac{\pi}{2}-o(1) =\angle \pmb{Z}_{j-1}\pmb{X}''\pmb{Z}_j < 
      \angle \pmb{Z}_{j-1}\pmb{X}_j\pmb{Z}_j  < \angle \pmb{Z}_{j-1}\pmb{X}'\pmb{Z}_j  
      <
      \angle \pmb{Z}_{j-1}\pmb{Z}_j\pmb{Z}'  = \frac{3\pi}{4}.
      $
       }

\end{figure}

\begin{lemma}\label{Lemma3}
For any affine plane $S$ of $\mathbb{R}^3$ such that $\pmb{Z}_{j-1}, \pmb{Z}_j\in S$,
 there exists $ \pmb{X}_j \in \mathcal{U}_{j}\cap S$ such that 
 \begin{equation}\label{AW}
 \mathcal{W}_j = \mathcal{W}_j(\pmb{X}_j) \subset \mathcal{U}_{j}
 \end{equation}
 and
  \begin{equation}\label{AW1}
 \mathcal{W}_j  \cap \mathcal{W}_j(\pmb{Z}_j) = \varnothing.
 \end{equation}
 Moreover  
 if we suppose that
 $j$ is large enough and  $\mathcal{W}_{j-1}$ satisfies \eqref{AW1} with $j$ replaced by $j-1$ and
   \begin{equation}\label{AW2}
   \pmb{Z}_j\in \mathcal{W}_{j-1},
 \end{equation}
 then
 for any $\pmb{X} \in  \mathcal{W}_j $
 the angle between    vectors
$\pmb{Z}_{j} - \pmb{X}$ and $\pmb{Z}_{j-1}- \pmb{X}$ is 
greater that $\frac{\pi}{4}$ and less than $\frac{3\pi}{4}$.
\end{lemma}

\begin{proof}
To get \eqref{AW} we simply need to take a unit vector $\pmb{e}$ orthogonal to  the unit vector 
  $\pmb{e}' = \frac{\pmb{Z}_{j-1}-\pmb{Z}_j}{|\pmb{Z}_{j-1}-\pmb{Z}_j|}$ and such that 
  $\pmb{Z}_j + \pmb{e}\in S$
   and take
  $$
  \pmb{X}_j = \frac{1}{2}\cdot q_j^{-u_j}  \pmb{e} +  \frac{1}{4}\cdot q_j^{-u_j}  \pmb{e}'.
  $$
  Then  \eqref{AW}  and \eqref{AW1} follow from triangle inequality. 
  The condition on the angle between $\pmb{Z}_{j} - \pmb{X}$ and $\pmb{Z}_{j-1}- \pmb{X}$
can be established as follows. From \eqref{AW2} and the growth conditions ({\bf iv},{\bf v}) we see that in the triangle
$\pmb{X}\pmb{Z}_j\pmb{Z}_{j-1}$ the length of the side $\pmb{Z}_j\pmb{Z}_{j-1}$ is much larger than the length of the 
side $\pmb{X}\pmb{Z}_j$ and so the angle between  $\pmb{Z}_{j} - \pmb{X}$ and $\pmb{Z}_{j-1}- \pmb{X}$
is in the interval between $ \frac{\pi}{2} -o(1)$ and $\frac{3\pi}{4}$ (see Fig. 1).
\end{proof}

We need further notation. Consider the affine subspace
 \begin{equation}\label{R31}
 \mathbb{R}_1^3  =\{ (x_0,x_1,x_2,x_3)\in \mathbb{R}^4:\, x_0 = 1\}.
 \end{equation}
 For a point $ \pmb{X}= (X_1,X_2,X_3)\in \mathbb{R}^3$ we define the point
 $ \overline{\pmb{X}}= (1,X_1,X_2,X_3)\in \mathbb{R}^3_1$. For any set 
 $\mathcal{A}\subset\mathbb{R}^3$ we define
 $$
 \overline{\mathcal{A} } = \{ \overline{\pmb{X}}:\,\, \pmb{X} \in \mathcal{A}\} \subset \mathbb{R}^3_1.
 $$
 Of course we can apply Lemma \ref{Lemma3} to the similar objects in affine subspace $\mathbb{R}^3_1$ instead of $\mathbb{R}^3$.
 
  The proof of Theorem \ref{Thm2} is based on an inductive construction. Starting from an integer  vector
  $$ \pmb{z}_1 = (q_1,a_{1,1},a_{1,2}, a_{1,3})\,\,\,\,\,\,\text{ with large }\,\,\,\,\,\,q_1
  $$
 we construct a sequence of integer vectors $\pmb{z}_\nu , \nu =1,2,3,...$ with \eqref{iv} such that  any collection \eqref{collek} of successive vectors satisfy conditions  ({\bf i}) - ({\bf vii}), neighborhoods form a nested sequence
 $$ 
  \mathcal{U}_1\supset \mathcal{W}_1\supset \cdots \supset  \mathcal{U}_{\nu-1}\supset\mathcal{W}_{\nu-1}\supset   \mathcal{U}_{\nu}\supset \mathcal{W}_\nu\supset  \cdots
  $$
  with non-empty intersection and for every $\pmb{\eta}\in  \mathcal{U}_{2\nu+3}$ the  collection
  \eqref{collek} will be a set of five successive best approximation vectors. Then for the unique point of the intersection $ \cap_\nu  \, \mathcal{U}_{\nu}$ the sequence  $\pmb{z}_\nu , \nu =1,2,3,...$ will be the sequence of {\it all best approximation vectors starting from} $\pmb{z}_1$. In our construction we are not interested in the best approximation vectors before $\pmb{z}_1$. Limit vector
   $\pmb{\alpha}$ will satisfy  the desired conditions. For our purpose it is important to carry out the construction with a given large value of $q_1$, and so we should take care on the constants in signs $\asymp$ and $\ll$ our inductive proof. We emphasis here that all the constants in these signs
 involved in our process do not depend neither on $q_1$ nor on the number $\nu$ of the inductive step.    
 
 To start the inductive process we construct three first vectors $\pmb{z}_1,\pmb{z}_2, \pmb{z}_3$ satisfying all the necessary conditions.
 This will be done in Section \ref{belov} below.
 Then we start an inductive step which for given 
 vectors $\pmb{z}_{2\nu+1},\pmb{z}_{2\nu+2}  , \pmb{z}_{2\nu+3}$ 
 constructs two next vectors
 $\pmb{z}_{2\nu+4}$ and $\pmb{z}_{2\nu+5} $. This will be done in Section \ref{8} below.
 Namely  we will explain there how for given collection 
  \eqref{collek} satisfying  ({\bf i}) - ({\bf vii}) and for given  neighborhood $\mathcal{U}_{2\nu+3}$ satisfying the inductive assumption
one can construct the collection for the next value of $\nu$ satisfying the similar properties and the next neighborhoods   satisfying all the necessary conditions. We will describe the construction of the
vectors $\pmb{z}_{2\nu+4}$ and $\pmb{z}_{2\nu+5} $  and to verify all the properties of these vectors and the corresponding neighborhoods
$\mathcal{U}_{2\nu+4}\supset \mathcal{U}_{2\nu+5}$. This procedure corresponds to a general concept of \emph{going down} and \emph{going up}.

\section{Base of induction}\label{belov}

Here we explain how to construct the first three vectors of our inductive procedure.
We need the standard Davenport-type argument (see \cite{CasselsGeometry}, Ch. I, \S 2.4) which we explain in the next subsection. Then in Subsection \ref{three} we describe the construction itself.

\subsection{Primitive vectors}

An integer vector $\pmb{e} = (e_1,e_2,e_3,e_4) \in \mathbb{Z}^4$ ic called {\it primitive} if ${\rm g.c.d.} (e_1,e_2,e_3,e_4) =1$.
We define a couple of integer vectors $\pmb{e},\pmb{f}\in \mathbb{Z}^4$ or a triple of integer vectors
$\pmb{e},\pmb{f}, \pmb{g} \in \mathbb{Z}^4$ to be {\it primitive}
if it can be completed by integer vectors to a basis of the whole lattice $\mathbb{Z}^4$, that is 
$$
\langle \pmb{e},\pmb{f}\rangle_\mathbb{Z} =  \langle \pmb{e},\pmb{f}\rangle_\mathbb{R}\cap \mathbb{Z}^4\,\,\,\,\,\,\,
\text{or}
\,\,\,\,\,\,\,
\langle \pmb{e},\pmb{f},\pmb{g} \rangle_\mathbb{Z} =  \langle \pmb{e},\pmb{f},\pmb{g}\rangle_\mathbb{R}\cap \mathbb{Z}^4.
$$

\begin{lemma}\label{LemmaX} 
There exists an absolute constant $c>0$ such that for any point $ \pmb{\beta} = (\beta_1,\beta_2,\beta_3,\beta_4) \in \mathbb{R}^4$  with
$ |\pmb{\beta}|\ge 3$
 there exists an integer point  $\pmb{b} = (b_1,b_2,b_3,b_4) \in  \mathbb{Z}^4$ such that
\begin{equation}\label{prim1}
 |\pmb{\beta} - \pmb{b}| \le c \cdot \frac{\log  |\pmb{\beta}|}{\log\log  |\pmb{\beta}|}
\end{equation}
and the triple
$$
\pmb{e}_1 = (1,0,0,0), \,\,\,\,\, \pmb{e}_2 = (0,1,0,0), \,\,\,\,\, \pmb{b} = (b_1,b_2,b_3,b_4)
$$
is primitive. Obviously, the couple $
\pmb{e}_1,   \pmb{b} 
$ will be primitive also.
 \end{lemma}
 
  \begin{proof}
  
  By Erd\H{o}s theorem  (see formula (3) from \cite{E} and discussion in \cite{MoPrim})  there exists  absolute constant $c>0$ such that there exists a primitive point  
  $$ (b_3^*, b_4^*) \in \mathbb{Z}^2, \,\,\,\,\,\,\, {\rm g.c.d}(b_3^*, b_4^*) = 1$$
   with 
  $$
  \max_{j=3,4} |\beta_j - b_j^* | <  \frac{c}{2} \cdot \frac{\log |\pmb{\beta}|}{\log\log |\pmb{\beta}|}.
  $$
  Put $ \pmb{b}^* = (0,0, b_3^*,b_4^*)$.
  Then the triple
  $$
\pmb{e}_1 ,\,\,\,\,\, \pmb{e}_2 ,\,\,\,\,\,    b_1\pmb{e}_1 +b_2 \pmb{e}_2 +\pmb{b}^*
$$
is primitive for all integers $b_1, b_2\in \mathbb{Z}$.
It is clear that we can choose $b_1, b_2$ to satisfy (\ref{prim1}).$\Box$

  \end{proof}

  \begin{lemma}\label{LemmaY}
  Let $\pmb{e} = (e_1,e_2,e_3,e_4)
  \in \mathbb{Z}^4$ be a primitive integer vector and 
   $$ \pmb{\beta} = (\beta_1,\beta_2,\beta_3,\beta_4),\,\,
    \pmb{\gamma} = (\gamma_1,\gamma_2,\gamma_3,\gamma_4) \in \mathbb{R}^4$$
    be such that for  their norms
    we have the bounds
    \begin{equation}\label{conditions}
 |\pmb{e} |\ge 3, \,\,\,\,\,      \frac{|\pmb{\beta}| \, \log\log |\pmb{\beta}| }  {\log |\pmb{\beta}| } \ge 200 c |\pmb{e} |, \,\,\,\,\,
      \frac{|\pmb{\gamma}|\, \log\log |\pmb{\gamma}|}  {\log |\pmb{\gamma}|} \ge 2000 c |\pmb{\beta}|,
         \end{equation}
         (here constant $c$ is defined in Lemma \ref{LemmaX}).
    Then there exist integer vectors
    $$
    \pmb{f} = (f_1,f_2,f_3,f_4) ,\,\, \pmb{g} = (g_1,g_2,g_3,g_4) \in \mathbb{Z}^4
    $$
    such that the triple $ \pmb{e},\pmb{f},\pmb{g}$ is primitive and
  \begin{equation}\label{ine0}
    |\pmb{\beta} - \pmb{f}| \le\frac{\frak{|\pmb{\beta}|}}{10},\,\,\,\,\,\,\,\,\,\,
    |\pmb{\gamma} - \pmb{g}| \le \frac{\frak{ |\pmb{\gamma}|}}{10}
    \end{equation}
  \end{lemma} 
  
    \begin{proof}
  Vector $\pmb{e}_1 = \pmb{e}$ can be completed by three vectors
    $\pmb{e}_j = (e_{1,j},e_{2,j},e_{3,j},e_{4,j})\in \mathbb{Z}^4$ to a basis of the whole lattice $\mathbb{Z}^4$  in such a way that
    \begin{equation}\label{ee}
    \max_{1\le i,j\le 4}  |e_{i,j}| \le 8|\pmb{e}|.
    \end{equation}
    Then
   $$
    \pmb{\beta} = \sum_{j=1}^4 \beta_j'\pmb{e}_j,\,\,\,\,
    \text{with}
    \,\,\,\,
    \beta_j' \in \mathbb{R},
    $$
    and   by (\ref{ee}) and Cramer theorem we get  the upper bound.
   
     \begin{equation}\label{Kramer}
    \max_{1\le j \le 4} |\beta_j'| \le \frak{b} = 200 |\pmb{\beta}|\cdot|\pmb{e}|^3.
    \end{equation}
    Now we use Lemma \ref{LemmaX} to complete the vector $ \pmb{e} = \pmb{e}_1$ to a primitive couple $\pmb{e},\pmb{f}$ with
    \begin{equation}\label{ine1}
    \pmb{f} =  (f_1,f_2,f_3,f_4) =\sum_{j=1}^4 b_j' \pmb{e}_j,\,\,\,\,\,
     |
    \beta_j'-b_j'| \le c \cdot \frac{\log \frak{b}}{\log\log \frak{b}}.
    \end{equation}
    But $   \pmb{f} - \pmb{\beta} =  \sum_{j=1}^4 (b_j' -\beta_j')\pmb{e}_j$ and so formulas (\ref{ee},\ref{Kramer},\ref{ine1})  and (\ref{conditions}) give us
   $$
       |\pmb{\beta}-\pmb{f}|\le 
       16|\pmb{e}|\cdot c\cdot \frac{\log \frak{b}}{\log\log \frak{b}}\le
       20c\cdot|\pmb{e}|\cdot \frac{\log |\pmb{\beta}|}{\log\log |\pmb{\beta}|} \le\frac{|\pmb{\beta}|}{10},
   $$ 
    and this is the
   first inequality from (\ref{ine0}).
   We should note that 
  $|\pmb{f} |
  \le 2|\pmb{\beta}| 
$.
   
   Then we consider primitive couple $ \pmb{e}_1^* = \pmb{e},\pmb{e}_2^* = \pmb{f}$ and complete it to a basis
   of the whole lattice $\mathbb{Z}^4$ by  two vectors
   $\pmb{e}_j 
   ^*= (e_{1,j},e_{2,j},e_{3,j},e_{4,j})\in \mathbb{Z}^4, \, j=3,4 $ with
    \begin{equation}\label{ee1}
    \max_{1\le i,j\le 4}  |e_{i,j}^*| \le 4   \max (|\pmb{e}|,|\pmb{f}|) \le 8|\pmb{\beta}|.
    \end{equation}
   Now
   $$
     \pmb{\gamma} = \sum_{j=1}^4 \gamma_j'\pmb{e}_j^*,\,\,\,\,
    \text{with}
    \,\,\,\,
 \gamma_j' \in \mathbb{R}
    $$
    and
        \begin{equation}\label{Kramer1}
    \max_{1\le j \le 4} |\gamma_j'| \le \frak{c} = 1000 \frak{|\pmb{\gamma}|}|\pmb{\beta}|^3.
    \end{equation}
    Now we use Lemma \ref{LemmaX}
    to complete the primitive couple $\pmb{e} = \pmb{e}_1^*, \pmb{f} =\pmb{e}_2^*$  by vector 
    \begin{equation}\label{ccc}
     \pmb{g} =  (g_1,g_2,g_3,g_4) =\sum_{j=1}^4 c_j' \pmb{e}_j^*,\,\,\,\,\,
    \max_{1\le j \le 4} |
    \gamma_j'-c_j'| \le c \cdot \frac{\log \frak{c}}{\log\log \frak{c}}.
    \end{equation}
    to a primitive triple  $\pmb{e}, \pmb{f} , \pmb{g}$. By (\ref{ee1},\ref{Kramer1},\ref{ccc})  and (\ref{conditions}) we deduce upper bound for  from (\ref{ine0})
    
      \end{proof}

\subsection{The initial triple $\pmb{z}_1,\pmb{z}_2,\pmb{z}_3$}\label{three}

We take arbitrary primitive integer vector $ \pmb{z}_1 = (q_1,a_{1,1}, a_{1,2}, a_{1,3})$ with $q_1 \ge 1$ and $ \max_{1\le j \le 3} |a_{1,j}| \le \frac{q_1}{2}$ and consider the corresponding rational point $\pmb{Z}_1 =\left( \frac{a_{1,1}}{q_1},  \frac{a_{1,2}}{q_1}, \frac{a_{1,3}}{q_1}\right)\in \mathbb{R}^3$. We take arbitrary point 
$\pmb{X}' =(X_1',X_2',X_3') \in \mathbb{R}^3$ such that the distance between $\pmb{Z}_1$  and  $\pmb{X}$ is 
\begin{equation}\label{dis1}
|\pmb{X}' - \pmb{Z}_1| = q_1^{-(1+\lambda\frak{g}_{1,0})}.
\end{equation}
Then 
we take a point
$\pmb{X}'' =(X_1'',X_2'',X_3') \in \mathbb{R}^3$ to satisfy the conditions
\begin{equation}\label{dis2}
|\pmb{X}'' - \pmb{X}'| = q_1^{-\frak{g}_{1,0}(1+\lambda\frak{g}_{1})}
\end{equation}
and
\begin{equation}\label{dis3}
\angle  \pmb{Z}_1 \pmb{X}''  \pmb{X}' = \frac{\pi}{2}.
\end{equation}
We want to apply Lemma \ref{LemmaY} for the vectors
$$\pmb{e} = \pmb{z}_1, \,\,\,\,\,\pmb{\beta} =  q_1^{\frak{g}_{1,0}}\overline{\pmb{X}}',
 \,\,\,\,\,
\pmb{\gamma} =  q_1^{\frak{g}_{1,0}\frak{g}_{1}} \overline{\pmb{X}}''.$$
It is clear that for large values of $q_1$ conditions (\ref{conditions}) will be satisfied.
So Lemma \ref{LemmaY} gives us a primitive triple
$$\pmb{z}_1=\pmb{e}, \pmb{z}_2=\pmb{f}, \pmb{z}_3 = \pmb{g}$$
 satisfying (\ref{ine0}). From the choice of $\pmb{\beta}$ and (\ref{ine0}) we see that 
 \begin{equation}\label{ppo1}
 \frac{q_1^{\frak{g}_{1,0}}}{2}\le q_2\le2 {q_1^{\frak{g}_{1,0}}},\,\,\,\,\,\,\,
  \frac{q_1^{\frak{g}_{1,0}\frak{g}_{1}}}{2}\le q_3\le2 {q_1^{\frak{g}_{1,0}\frak{g}_{1}}}.
 \end{equation}
  
\vskip+0.3cm

First of all, in the construction described we can take $q_1$ to be arbitrary large.
Then we check all the conditions from ({\bf i}) - ({\bf vii}) which deal just with three vectors $\pmb{z}_1,\pmb{z}_2,\pmb{z}_3$.
Let us mention them below.

 1) Conditions  ({\bf ii}) and ({\bf iii})  are not applicable to to triple $\pmb{z}_1,\pmb{z}_2,\pmb{z}_3$.

2) Condition ({\bf i}) is satisfied as  vectors $\pmb{z}_1,\pmb{z}_2,\pmb{z}_3$ form a primitive triple.

3) Conditions ({\bf iv}) and  ({\bf v}) turn into  $ q_3 \asymp q_2^{\frak{g}_1}$,$\zeta_2 = q_2 |\pmb{Z}_2-\pmb{Z}_3|\asymp q_3^{-\lambda} $
and.
They follow from (\ref{ppo1}) for the denominators, and from (\ref{dis1},\ref{dis2}) and (\ref{ine0}) for the approximation.

4) Condition ({\bf vi}) follows from (\ref{dis3}) and the inequalities (\ref{ine0}) for the approximation.

5) As for the condition  ({\bf vii}) which gives upper and lower  bounds (\ref{wowo}) for the fundamental volume of the lattice
$\langle \pmb{z}_1,\pmb{z}_2,\pmb{z}_3\rangle_\mathbb{Z}$,
the upper bound of the type (\ref{tw1}) while the lower bound immediately follows from the condition ({\bf vi}) on the angle and Lemma \ref{Lemma6} below.

 \begin{lemma} \label{Lemma6}  Consider the three-dimensional parallelepiped 
  $$
  \Pi = \{ \pmb{z} = \lambda_1 \frak{z}_1+ \lambda_2  \frak{z}_2 + \lambda_3  \frak{z}_3,\,\,\,
  0\le \lambda_j \le 1 \}
  $$
  generated by  linearly independent vectors $ \frak{z}_j = (q_j,a_{1,j}, a_{2,j},a_{3,j}) 
   \in \mathbb{R}^4 , \, j =1,2,3$.
  Let $\frak{Z}_j = \left(\frac{a_{1,j}}{q_1},\frac{a_{2,j}}{q_2}, \frac{a_{3,}j}{q_3}\right)
  $ be the corresponding points on $\mathbb{R}^3$.
  Then
  $$
  {\rm Vol} \, \Pi \gg q_1q_2q_3\cdot |\frak{Z}_1- \frak{Z}_3|\cdot |\frak{Z}_2-\frak{Z}_3|\cdot \sin \varphi
  $$
  where $\varphi $ is the angle between $\frak{Z}_1- \frak{Z}_3$ and $\frak{Z}_2- \frak{Z}_3$.
  \end{lemma}
  
\begin{proof}We consider linear transformation $\mathcal{G}$  such that 
$$
\mathcal{G}\frak{z}_j =  \overline{\frak{Z}}_j = \left(1, \frac{a_{1,j}}{q_1},\frac{a_{2,j}}{q_2}, \frac{a_{3,}j}{q_3}\right),
\,\,\,\,
j = 1,2,3.
$$
Then  
$$
 {\rm Vol} \, \mathcal{G} \Pi = \frac{1}{q_1q_2q_3} {\rm Vol} \, \Pi .
 $$
 But  $\mathcal{G} \Pi$ contains simplex with vertices $\pmb{0},  \overline{\frak{Z}}_1, \overline{\frak{Z}}_2, \overline{\frak{Z}}_3$ whose volume is greater than
 $$ \frac{ |\frak{Z}_1- \frak{Z}_3|\cdot |\frak{Z}_2-\frak{Z}_3|\cdot \sin \varphi}{6}.$$
 Lemma is proven.
 \end{proof}

To conclude this subsection we should outline that 
all the constants in signs $\asymp$ in conclusion 3)  do not depend on $q_1$ because the constants in Lemma \ref{LemmaY} are absolute.
Also for the conclusion 5) and  the constants in the condition   ({\bf vii}), they also do not depend on $q_1$, because (\ref{tw1}) deals with  absolute constants and the constants in Lemma \ref{Lemma6} are also absolute.

\section{Two- and three-dimensional lattices}\label{7}

We consider the two-dimensional subspace
$$
 L = \langle\pmb{z}_{2\nu+2},\pmb{z}_{2\nu+3}
\rangle_\mathbb{R}  
$$
the two-dimensional lattice
$$
 \Lambda = \langle\pmb{z}_{2\nu+2},\pmb{z}_{2\nu+3}
\rangle_\mathbb{Z}  = \mathbb{Z}^4 \cap L\subset L
$$
as well as the three-dimensional subspace
$$
G=  \langle\pmb{z}_{2\nu+1},\pmb{z}_{2\nu+2},\pmb{z}_{2\nu+3}
\rangle_\mathbb{R}
$$
and  the lattice
$
\Gamma $ defined in ({\bf vii}).
Let $ d$ be the two-dimensional fundamental volume of lattice $\Lambda$ and $\Delta$ be  the three-dimensional fundamental volume of lattice $\Gamma$ defined in condition ({\bf vii}). Then
\begin{equation}\label{dD}
d\asymp \zeta_{2\nu+2}q_{2\nu+3}\asymp q_{2\nu+3}^{1-\lambda},\,\,\,\,\,\,
\Delta \asymp \zeta_{2\nu +1}\zeta_{2\nu+2} q_{2\nu +3} \asymp q_{2\nu+2}^{-\lambda}q_{2\nu+3}^{1-\lambda}
\asymp q_{2\nu+3}^{1-\lambda-\frac{\lambda}{\frak{g}_1}},\,\,\,\,\,\,
\frac{\Delta}{d} \asymp q_{2\nu+2}^{-\lambda}.
\end{equation}

We need to consider the affine plane
 $$
 L_{1}
=  L+\pmb{z}_{2\nu+1} \subset G
$$
and the affine lattice
 $$
 \Lambda_{1}
=  \Lambda+\pmb{z}_{2\nu+1} \subset \Gamma.
$$
We should note that the Euclidean distance between the parallel planes
$L$ and $L_1$ is equal to
$\frac{\Delta}{d}.
$

Also we should deal with three-dimensional affine subspace
 $$
 G_{1}
=  G+\pmb{z}_{2\nu}  \subset \mathbb{R}^4
$$
and the affine lattice
 $$
 \Gamma_{1}
=  \Gamma+\pmb{z}_{2\nu} \subset \mathbb{Z}^4.
$$
The Euclidean distance between  $G$ and $G_1$ is equal to
$\frac{1}{\Delta}.
$

The parallelogram
$$
\mathcal{P} = \{ \pmb{z} =(x,y_1,y_2,y_3) = (x,\pmb{y})\in L :
0\le x \le q_{2\nu+3},\,\, |\pmb{y} - x\pmb{Z}_{2\nu+3}|\le \zeta_{2\nu+2}\}
\subset L
$$
contains  integer points $ \pmb{0},\pmb{z}_{2\nu+2},\pmb{z}_{2\nu+3},  \pmb{z}_{2\nu+3}-\pmb{z}_{2\nu+2}$
and hence a fundamental domain of the lattice $\Lambda$,  while the cylinder 
$$
\mathcal{C}= \{ \pmb{z} =(x,y_1,y_2,y_3)\textcolor{blue}{=} (x,\pmb{y})\in G :
0\le x \le q_{2\nu+1}+q_{2\nu +2}+q_{2\nu+3}, \,\,
|\pmb{y}-x\pmb{Z}_{2\nu+3}|
\le 3\zeta_{2\nu+1}\} \supset
$$
$$\supset
\{ \pmb{z} =(x,y_1,y_2,y_3)= \lambda_1 \pmb{z}_{2\nu+1}+ \lambda_2 \pmb{z}_{2\nu+2}
+ \lambda_3 \pmb{z}_{2\nu+3}, \,\, 0\le \lambda_1,\lambda_2,\lambda_3\le 1\}
$$
contains a fundamental domain of the lattice $\Gamma$.

 \section{Step of induction}\label{8}

 In the proof of Theorem \ref{Thm2} for case $k=1$, conditions 
({\bf vi},{\bf v}) and the exponents $\frak{g}_1$ and $\frak{g}_{1,0}$ defined in Lemma \ref{Lemma1} and \eqref{g11} play a crucial role.

 Now we assume that the vectors $\pmb{z}_{2\nu+1},\pmb{z}_{2\nu+2}  , \pmb{z}_{2\nu+3}$  are constructed.  We should construct vector   $\pmb{z}_{2\nu+4}$ ({\bf Stage 1}) and
 vector  $\pmb{z}_{2\nu+5}$ ({\bf Stage 2}).
  Moreower, for the collection  \eqref{collek} we assume that it satisfies ({\bf i}) - ({\bf vii}); also we suppose that for every  $\pmb{\xi} \in \mathcal{U}_{2\nu+3}$
  the collection \eqref{collek} 
 contains of just successive five best approximation to $\pmb{\xi}$.  
 
\subsection{Stage 1}

 First, we construct integer vector $\pmb{z}_{2\nu+4}$. 
 For  $\pmb{z}_{2\nu+3}$ Lemma 3
 applied in $\mathbb{R}_1^3$ with $ S = \mathbb{R}_1^3\cap G$
  gives us the point $X_{2\nu+3}$ and the neighborhood $\mathcal{W}_{2\nu+3}$.

  Then we consider the central projection $\frak{P}_1$  with center ${\pmb{0}}$ onto the affine subspace $L_1$, that is for $\pmb{x}\in \mathbb{R}_1^3$ its image (if it exists)  is defined as intersection of the line ${\rm span} \, \pmb{x}$ with subspace $L_1$. So $\frak{P}_1$    is  a map from $S\setminus L$ to $L_1$.

  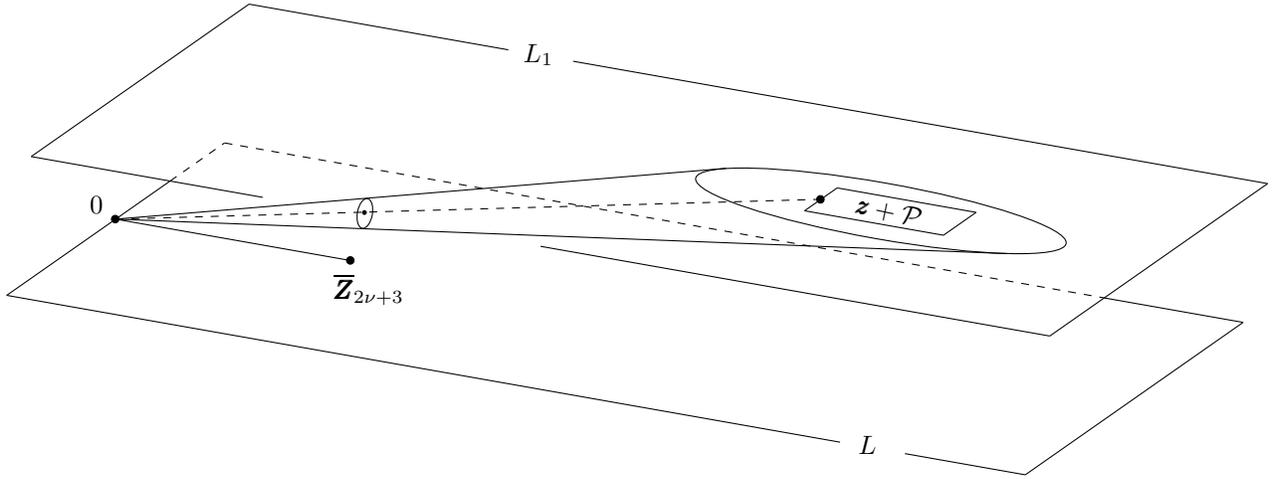
\begin{figure}[h]
  \centering

\begin{tikzpicture}[scale=1.25]

\draw[rotate=-10] (-2.2,2) -- (-5,2);
\draw[rotate=-10] (-1.5,2) -- (6,2);

\draw[rotate=-10] (-7,0) -- (-4.5,0);
\draw[rotate=-10] (-1.5,0) -- (4,0);

\draw[rotate=-10] (-5,2) -- (-7,0);
\draw[rotate=-10] (4,0) -- (6,2);

    \draw (2,1) [rotate=-10] ellipse (2cm and 0.3cm);
    
       \draw (-3.4,0.02) [rotate=-10] ellipse (0.08cm and 0.16cm); 
    
    \draw[rotate=-10,dashed] (-5,0.5) -- (4.5,0.5);
     \draw[rotate=-10] (6,0.5) -- (4.5,0.5);

\draw[rotate=-10] (-7,-1.5) -- (2,-1.5);

\draw[rotate=-10] (4,-1.5) -- (2.7,-1.5);

\draw[rotate=-10] (-5.5,0) -- (-7,-1.5);
\draw[rotate=-10] (4,-1.5) -- (6,0.5);

   \draw[rotate=-10,dashed] (-5.5,0) -- (-5,0.5);
\draw[rotate=-10] (-6,-0.5) -- (-3.5,-0.5);

\draw[rotate=-10] (-6,-0.5) -- (0.3,1.16);
\draw[rotate=-10] (-6,-0.5) -- (3.4,0.785);

 \node[draw=black,fill=black,circle,inner sep=1pt] at (-6, 0.55) {};
 
 \node[draw=black,fill=black,circle,inner sep=1pt] at (-3.5, 0.11) {};
 
 \node[draw=black,fill=black,circle,inner sep=1pt] at (1.5, 0.76) {};
 
 \draw[rotate=-10] (1.2,0.86) -- (1.5,1.16);
 
 \draw[rotate=-10] (2.7,0.86) -- (3,1.16);
 
 \draw[rotate=-10]  (1.5,1.16) -- (3,1.16);
 
  \draw[rotate=-10] (1.2,0.86) -- (2.7,0.86);
 
  %  \draw[rotate=-10] (0.77,0.77) -- (4,1);
  
   % \draw[rotate=-10] (4,1) -- (1.28,1.28);
   
   % \draw[rotate=-10] (0.86,0.86) -- (2.05,0.86);
   
      %\draw[rotate=-10] (1.17,1.17) -- (2.33,1.17);
      
      % \draw[rotate=-10]  (2.05,0.86) -- (2.33,1.17);
       
         \draw[rotate=-10,dashed] (-6,-0.5) -- (1.3,1);
          
           \node[draw=black,fill=black,circle,inner sep=0.5pt] at (-3.35, 0.62) {};

 \draw (-6.2, 0.7)   node { $0$};
 
  \draw (2, -1.85)   node { $L$};
  
    \draw (-1.5, 2.3)   node { $L_1$};
    
       \draw (-3.3, -0.2)   node { $\overline{\pmb{Z}}_{2\nu+3}$};

            \draw(2.23, 0.63) node  [rotate=-12]   { $\pmb{z} + \mathcal{P}$} ;

\end{tikzpicture}

      \caption{To Lemma 4: small ellipse at the left side close to the point $\overline{\pmb{Z}}_{2\nu+3}$  is the neighborhood $\overline{\mathcal{W}}_{2\nu+3}\subset \mathbb{R}_1^3$  with center $\overline{\pmb{X}}_{2\nu+3}$ and  the large ellipse at the right side is the image  $\frak{P}_1(\overline{\mathcal{W}}_{2\nu+3})  \subset L_1$.}

\end{figure}
  
 \begin{lemma}\label{Lemma4}. Suppose that $q_1$ is large enough. Then the image $ \frak{P}_1( \overline{\mathcal{W}}_{2\nu+3 })$ contains a shift of  a fundamental  domain of lattice $\Lambda$.
 \end{lemma}
 
 \begin{proof}
 The geometry of the proof is illustrated on Fig. 2.\\
   
Consider the unit vector
 $\pmb{e}_0=
 \frac{\pmb{z}_{2\nu+3}}{|\pmb{z}_{2\nu+3}|}
 $
 and the unit vector
 $\pmb{e}_1\in   G, |\pmb{e}_1| = 1$ which is orthogonal to $L$
 such that    $\frac{\Delta}{d}\pmb{e}_1\in L_1$. We also need a unit vector $ \pmb{e}_2 = (0,u,v,w)\in L$ which belongs to subspace $L$,  and  such that both vectors $\pmb{e}_2$ and $\pmb{z}_{2\nu+2}$ lie in the same half-subspace of $L$ with respect to the one-dimensional subspace $ \langle\pmb{z}_{2\nu+3}\rangle_\mathbb{R} \subset L$. Unit vectors $\pmb{e}_j, j =0,1,2$ are defined uniquely by the conditions above.  
    
    As $ \overline{\pmb{X}}_{2\nu+3} \in S\subset G$ the image
    $\pmb{z} =\frak{P} (\overline{\pmb{X}}_{2\nu+3})\in L_1$ is well defined and satisfies
    \begin{equation}\label{zee}
    \pmb{z} =  r \pmb{e}_0+
 s\pmb{e}_1+t \pmb{e}_2.
    \end{equation}
Calculations show that 
\begin{equation}\label{show}
 r 
 \asymp  \frac{\Delta}{d}\cdot q_{2\nu+3}^{1+\lambda \frak{g}_{1,0}} 
 ,\,\,\,\,\,\,
 s =\frac{\Delta}{d},
 \,\,\,\,\,\,
 t\asymp  \frac{\Delta}{d},
\end{equation}
and the constant in  sign $\asymp$ are absolute and do not depend on the values of $q_j$. We will writhe down explicit constants in the exposition below.

Indeed, the equality for $s$ from (\ref{show})  follows from the fact that the unit vector $\pmb{e}_2\in G$ is orthogonal to $L$ and the distance between the subspaces $L$ and $L_1$ is equal to $ \frac{\Delta}{d}$. By the way, the distance between the subspaces $L$ and $L_1$ is just the distance from the point $\pmb{z}$ to the subspace $L$.
Now we define parameters $ r_*,s_* t_*$ by
$$
\overline{\pmb{X}}_{2\nu+3}  =  r_* \pmb{e}_0+
 s_*\pmb{e}_1+t_* \pmb{e}_2.
 $$
So, in particular $s_*$ is the distance from the point $\overline{\pmb{X}}_{2\nu+3} $ to the subspace $L$. Then
\begin{equation}\label{show1}
\frac{r}{r_*} = \frac{s}{s_*}=\frac{t}{t_*}.
\end{equation}
But by (\ref{fri}) we see that 
\begin{equation}\label{show2}
q_{2\nu+3} \le r_*= |\pmb{z}_{2\nu+3}|\le 4 q_{2\nu+3}.
\end{equation}
Now we consider the point of the form 
$$
\overline{\pmb{X}}_{**} = (1, X^1,X^2, X^3) \in \mathbb{R}^3_1 \cap G
$$
(recall that $\mathbb{R}^3_1$ us defined in (\ref{R31}))
which is the orthogonal projection of the point $ \overline{\pmb{X}}_{2\nu+3} \in \mathbb{R}^3_1\cap G$ onto the line
$  \mathbb{R}^3_1\cap L\subset  \mathbb{R}^3_1\cap G$. We consider triangle $ \overline{\pmb{X}}_{2\nu+3}\overline{\pmb{Z}}_{2\nu+3} \overline{\pmb{X}}_{**} \subset  \mathbb{R}^3_1\cap G$.

Let $ t_{**}$ be the distance between points $ \overline{\pmb{X}}_{2\nu+3}$ and $
\overline{\pmb{X}}_{**} $. By taking into account the angle 
 $\varphi_*$ between vectors 
$ \pmb{z}_{2\nu+3} = (q_{2\nu+3}, a_{1,2\nu+3}, a_{2,2\nu+3}, a_{3,2\nu+3})$ and $(1,0,0,0)$, we see that 
\begin{equation}\label{show5}
1\ge \frac{t_*}{t_{**}} \ge \cos \varphi_*
= \frac{q_{2\nu+3}}{ |\pmb{z}_{2\nu+3}|}\ge \frac{1}{4}.
\end{equation}
Because of $  \overline{\pmb{X}}_{2\nu+3}\in  \overline{\mathcal{W}}_{2\nu+3 }\subset \overline{\mathcal{U}}_{2\nu+3 }$
 we conclude (see Fig. 1)  that 
 \begin{equation}\label{show6}
 \frac{q_{2\nu+3}^{-u_{2\nu+3}}}{10}
  \le t_{**}\le q_{2\nu+3}^{-u_{2\nu+3}}.
\end{equation}
Let $ s_{**}$ be the distance between points $ \overline{\pmb{Z}}_{2\nu+3}$ and $
\overline{\pmb{X}}_{**} $. 
  We see that 
\begin{equation}\label{show3}
1\ge \frac{s_*}{s_{**}} \ge \cos \varphi_*
 = \frac{q_{2\nu+3}}{ |\pmb{z}_{2\nu+3}|}\ge \frac{1}{4},
\end{equation}
by (\ref{fri}).
 Again, because of $  \overline{\pmb{X}}_{2\nu+3}\in  \overline{\mathcal{W}}_{2\nu+3 }\subset \overline{\mathcal{U}}_{2\nu+3 }$
 we have
 \begin{equation}\label{show4}
 \frac{q_{2\nu+3}^{-u_{2\nu+3}}}{20}\le
  \left(\left( \frac{1}{10}{\sin(\pi/8)} -\frac{1}{10}\right) \cdot \cos (\pi/8) -\frac{1}{10}\right)q_{2\nu+3}^{-u_{2\nu+3}}\le s_{**}\le q_{2\nu+3}^{-u_{2\nu+3}}.
\end{equation}
Now (\ref{show1}) gives us $ r = r_* \frac{s}{s_*} = \frac{\Delta}{d} \cdot \frac{r_*}{s_*}$ and inequalites
(\ref{show2},\ref{show3},\ref{show4}) together with the definition (\ref{uji}) lead to the desired bound
$$
\frac{\Delta}{d}  \cdot q_{2\nu+3}^{1+\lambda \frak{g}_{1,0}} 
\le r \le 320 \cdot \frac{\Delta}{d}  \cdot q_{2\nu+3}^{1+\lambda \frak{g}_{1,0}} .
$$
Analogously, for $t$ from (\ref{show1},\ref{show3},\ref{show4},\ref{show5},\ref{show6})  we deduce inequalities
$$
\frac{1}{40}\cdot \frac{\Delta}{d} \le t \le 80\cdot \frac{\Delta}{d} .
$$

 From \eqref{dD} and  ({\bf iv})  and the definition \eqref{g11}  of $\frak{g}_{1,0}$ we see that 
  \begin{equation}\label{rst1}
  r\asymp  q_{2\nu+3}^{\frak{g}_{1,0}},
\,\,\,\,\,\text{
   meanwhile}\,\,\,\,\,
  t  \asymp q_{2\nu+2}^{-\lambda},
   \end{equation}
   with explicit absolute constants in signs $\asymp$.
  The image $ \frak{E}= \frak{P}_1 ( \overline{\mathcal{W}}_{2\nu+3 }) $ is an ellipse.
  The angle between the vector $\pmb{z}_{2\nu+3}$  and the 
 large major axis  of this ellipse is $O(q_{2\nu+3}^{-u_{2\nu+3}})$. The length of 
 the 
 large major axis is 
$\asymp r$. The section $ \{ x_0 =r\}$ of the ellipse  has length $\asymp t$. 
 From \eqref{rst1} we see (as $q_{2\nu+3}>q_{2\nu+2}>q_1$, meanwhile $q_1$  is large enough),
 that the large axis  fo $\frak{E}$ is $ \ge 4q_{2\nu+3}$ and the section  $ \{ x_0 =r\}$  has length $ \ge 4\xi_{2\nu+2}$. Calculations show that the vertices of the parallelogram
 $ \pmb{z}+ \mathcal{P}$ belong to ellipse $\frak{E}$.
Then, by convexity
 $$
 \pmb{z}+ \mathcal{P}\subset  \frak{E} = \frak{P}_1 ( \overline{\mathcal{W}}_{2\nu+3 }).
 $$
 Lemma is proven.
 \end{proof}
 
 Recall that $\Lambda$ is a two-dimensional lattice in the plane $L$.
 Now we consider the lattice $\Lambda_1=\Lambda +\pmb{z}_{2\nu+1}$ which belongs the the plane
 $L_1$
By Lemma 4,  $ \frak{P}_1 ( \overline{\mathcal{W}}_{2\nu+3 })$ contains a  shift of a fundamental domain of lattice $\Lambda$ and lattice $\Lambda_1$ is congruent to $\Lambda$. So the image 
  $ \frak{P}_1 ( \overline{\mathcal{W}_{2\nu+3 }})$ contains an integer point. This is just the point $\pmb{z}_{2\nu+4} = (q_{2\nu+4} , a_{1,2\nu+4}, a_{2,2\nu+4}a_{3,2\nu+4})$ which we are constructing. 
  
  So we have defined the next point $\pmb{z}_{2\nu+4}$.
  
 First of all we should note that, as $ \pmb{z}_{2\nu+4}  \in \Lambda_1$, it is of the form
 \begin{equation}\label{newform}
  \pmb{z}_{2\nu+4}  =  \pmb{z}_{2\nu+1}  + a_{\nu+1}\pmb{z}_{2\nu+2}  +b_{\nu+1}\pmb{z}_{2\nu+3} 
 \end{equation}
 with integers 
 $ a_{\nu+1},b_{\nu+1}$
 and we have ({\bf ii})  for the next step. Then
 $$
 \pmb{z}_{2\nu+4} \in \pmb{z}+\cal{P},
 $$
 where $\pmb{z}$ is of the form  (\ref{zee}) and coefficients $r,s,t$ satisfy (\ref{show}).
 Moreover for any point of parallelogram $\cal{P}$ its first coordinate is $O(q_{2\nu+3})$.
 So 
 from \eqref{zee} and   \eqref{rst1} we see that 
 $$
 q_{2\nu+4}\asymp  r \asymp q_{2\nu+3}^{\frak{g}_{1,0}}
 $$
 and this is just the new (second) inequality from the condition ({\bf v}) for the next step. 
 
 The condition $ \zeta_{2\nu+3} \asymp   q_{2\nu+4}^{-\lambda}$  follows from the inclusion
  $ \pmb{Z}_{2\nu+3}\in \mathcal{W}_{2\nu+3}  \subset \mathcal{U}_{2\nu+3}$
 and our choice of parameters.

 Note that  the approximation $\zeta_{2\nu+4}$ is still not defined. It will be defined at the next {\bf Stage 2} when we define the next vector $\pmb{z}_{2\nu+5}$.
 
  It is clear from triangle inequality  and  the asymptotics for parameters \eqref{rst1} that 
 $\mathcal{U}_{2\nu+4}\subset\mathcal{W}_{2\nu+3} $
 
 Now we must explain why this point 
  will be the next to $\pmb{z}_{2\nu+3}$ best approximation for all  points from the next neighborhood $\mathcal{U}_{2\nu+4}$ that is, for every $\pmb{\eta}\in \mathcal{U}_{2\nu+4}$  vectors $\pmb{z}_{2\nu+3}$ and $\pmb{z}_{2\nu+4}$ are successive best approximation vectors.

 First of all we show that the best approximation vectors for $\pmb{\eta}$ with 
  $ q_{2\nu+3}\le q\le q_{2\nu+4}$ should belong to the three-dimensional subspace $G$.
 Indeed,
 the whole lattice $\mathbb{Z}^4$ splits into three-dimensional affine sublattices
 $$
 \Gamma_k = \Gamma+ k \pmb{z}_{2\nu}\subset G_k = G+ k \pmb{z}_{2\nu}
 $$
  parallel to $\Gamma \subset
 G$. The distance between two parallel three-dimensional affine subspaces $ G_k$ and 
 $G_{k+1}$ is  equal to 
  $$
 \frac{1}{\Delta}\asymp \zeta_{2\nu} \asymp q_{2\nu+1}^{-\lambda},
$$
So, for any integer vector  
  $ \pmb{w} = (q,\pmb{a}) =(q,a_1,a_2,a_3)\in \mathbb{Z}^4 \setminus G$ one has
  $$ \min_{\pmb{w}'\in G} |\pmb{w}-\pmb{w}'| \ge \frac{1}{\Delta} .
  $$
  Now  for    $ \pmb{w} =(q,a_1,a_2,a_3)\in \mathbb{Z}^4 \setminus G $ with $ q_{2\nu+3}\le q\le q_{2\nu+4}$  and for any $\pmb{\eta}\in \mathcal{U}_{2\nu+3}$ we get
 $$
 |q\pmb{\eta} - \pmb{a}|
 \ge
      \min_{\pmb{w}'\in G} |\pmb{w}-\pmb{w}'|  - q\times  \textrm{diam}(\mathcal{U}_{2\nu+3}) =
   \frac{1}{\Delta}  - q\times  \textrm{diam}(\mathcal{U}_{2\nu+3}), 
 $$
 where $ \textrm{diam}(\mathcal{U}_{2\nu+3})$ is the diameter of the neighborhood $\mathcal{U}_{2\nu+3}$. But
 \begin{equation}\label{previous}
 q\times \textrm{diam}(\mathcal{U}_{2\nu+3}) \ll  q_{2\nu+4} q_{2\nu+3}^{-1-\lambda \frak{g}_{1,0}} \asymp
 q_{2\nu+3}^{-\lambda/\frak{g}_1} \asymp q_{2\nu+2}^{-\lambda} = o(q_{2\nu+1}^{-\lambda}),
 \end{equation}
 and so
  $$
     |q\pmb{\eta} -  \textcolor{blue}{ \pmb{a}}|> |q_{2\nu+3}\pmb{\eta} -\pmb{a}_{2\nu+3}|.
  $$
  In the last formula we establish the sign $>$ but not $\gg$. The explanation is that in the previous estimate (\ref{previous}) the constants in the signs $\asymp, \ll, o(\cdot)$ do not depend on $q_1$, meanwhile $ q_{2\nu+1} \ge q_1$ and we take $q_1$ to be large enough.
 Meanwhile in $G$ the only candidates to be new best approximation  vectors  in the range $ q_{2\nu+3}\le q\le q_{2\nu+4}$ are $\pmb{z}_{2\nu+4}$ and $\pmb{z}_{2\nu+4}-\pmb{z}_{2\nu+3}$
 (there is no integer points in $G$ in the strip between $L$ and $L_1$ and the distances form vectors  $\pmb{z}_{2\nu+3}$ and $\pmb{z}_{2\nu+4}-\pmb{z}_{2\nu+3}$ to the linear subspace generated by $\pmb{z}_{2\nu+4}$ are equal).
 Vector $\pmb{z}_{2\nu+4}$ is a best approximation vector indeed, and $\pmb{z}_{2\nu+4}-\pmb{z}_{2\nu+3}$ is not a best approximation vector as $\pmb{\eta}\in \mathcal{U}_{2\nu+4}$ and this angular neighborhood just separates  $\pmb{\eta}$ from
  $\pmb{z}_{2\nu+4}-\pmb{z}_{2\nu+3}$. So the condition $\pmb{\eta}\in \mathcal{U}_{2\nu+4}$  ensures that the vector $\pmb{z}_{2\nu+3}$ is closer than $\pmb{z}_{2\nu+4}-\pmb{z}_{2\nu+3}$  to the line generated by $\pmb{\eta}$, meanwhile $  q_{2\nu+4}-q_{2\nu+3}> q_{2\nu+3}$.
 
  Indeed, the points $\pmb{0}, \pmb{z}_{2\nu+2}, \pmb{z}_{2\nu+4}, \pmb{z}_{2\nu+4}-\pmb{z}_{2\nu+3}$
  form a parallelogram and so
  $$
  | q_{2\nu+3}{\overline{\pmb{Z}}_{2\nu+4}} -\pmb{z}_{2\nu+3}|=
    | q_{2\nu+3}{\overline{\pmb{Z}}_{2\nu+4}} -(\pmb{z}_{2\nu+4}-\pmb{z}_{2\nu+3})|.
  $$
  But
  we take $\pmb{\eta}\in \mathcal{U}_{2\nu+4}\subset \mathcal{W}_{2\nu+3}$ (see Fig. 1) and so 
  $\pmb{\eta}$ is "closer" to $\pmb{z}_{2\nu+3}$ than to $\pmb{z}_{2\nu+4} - \pmb{z}_{2\nu+3}$ and
  $$
  |q_{2\nu+3}\pmb{\eta} - \pmb{a}_{2\nu+3}|< 
    |(q_{2\nu+4} - q_{2\nu+3})\pmb{\eta} - (\pmb{a}_{2\nu+4}-\pmb{a}_{2\nu+3})|
.
  $$
  
  We have established all the approximation properties for all the vectors $\pmb{\eta}\in \mathcal{U}_{2\nu+4}$.
  
\subsection{Stage 2}

 Now we should construct vector $\pmb{z}_{2\nu+5}$ and  establish all the necessary properties of points 
  $\pmb{\eta}$ in the new neighbourhood $\mathcal{U}_{2\nu+5}$.

      First, we should note that by ({\bf{i}}) for the second triple and \eqref{newform}, the triple
  $\pmb{z}_{2\nu+2},\pmb{z}_{2\nu+3}, \pmb{z}_{2\nu+4}
 $
 form a basis of the lattice $\Gamma$, while the quadruple  
 $\pmb{z}_{2\nu},\pmb{z}_{2\nu+2},\pmb{z}_{2\nu+3}, \pmb{z}_{2\nu+4}
 $
 is a basis of $\mathbb{Z}^4$ (we take into account ({\bf iii})).
 
 We should note that {$\pmb{z}_{2\nu}\not \in G$}.
 Together with three-dimensional linear subspace $G$ we consider affine subspace
 $ G_1 = G+ {\pmb{z}_{2\nu}}$.
 The Euclidean distance between parallel subspaces $G$ and $G_1$ is equal to $\frac{1}{\Delta}$.
 We want to construct a new integer point $\pmb{z}_{2\nu+5}$  which belongs to the rational affine subspace $
 G_1$. The subspace $G_1$ contains a three-dimensional integer lattice  $G_1\cap\mathbb{Z}^4$. This gives us certain freedom to choose the point $\pmb{z}_{2\nu+5}$. However it would be enough  to localise the new point $\pmb{z}_{2\nu+5}$ 
close to a three-dimensional linear subspace $S$ 
 defined below.

We should apply Lemma \ref{Lemma3}. We use it for $ j = 2\nu +4$ and for points $\pmb{Z}_{\nu+4},\pmb{Z}_{\nu+3}$. We choose $S$ as follows.
    Let $G'\subset \mathbb{R}^4$ be three dimensional subspace with basis
  $ \pmb{z}_{2\nu+3}, \pmb{z}_{2\nu+4}, \pmb{f}_1$,
  where   $\pmb{f}_1, |\pmb{f}_1| = 1$  is the unit which is orthogonal to $G$ and
 such that    $\frac{1}{\Delta}\pmb{f}_1\in G_1$. 
 We define $ S = G'\cap \mathbb{R}_1^3$. Then Lemma \ref{Lemma3} gives a point $\pmb{X}_{2\nu+4}$ and a neighborhood $\mathcal{W}_{2\nu+4}$.
  
  We consider the central projection $\frak{P}_2$ with center $\pmb{0}$ onto the affine subspace $G_1$. So
   $\frak{P}_2 $ is a map from
     $\mathbb{R}^4\setminus  G$ to $G_1$.
 
 \begin{lemma} \label{Lemma5}
 Suppose that $q_1$ is large enough.
  Then the  image $ \frak{P}_2( \overline{\mathcal{W}}_{2\nu+4 })$ contains a shift of  a fundamental  domain of lattice $\Gamma$.
 \end{lemma}
 
 \begin{proof}
 We consider the unit vector
 $\pmb{f}_0=
 \frac{\pmb{z}_{2\nu+4}}{|\pmb{z}_{2\nu+4}|}
 $
 and the unit vector
 $\pmb{f}_1, |\pmb{f}_1| = 1$ which is orthogonal to $G$
defined above. 
 Again, we need a unit vector $ \pmb{f}_2 = (0,u,v,w)\in  \langle \pmb{z}_{2\nu+3}, \pmb{z}_{2\nu+4}  \rangle_{\mathbb{R}}$  such that both vectors $\pmb{f}_2$ and $\pmb{z}_{2\nu+3}$ lie in the same half-subspace of the two-dimensional subspace $\langle \pmb{z}_{2\nu+3}, \pmb{z}_{2\nu+4}  \rangle_{\mathbb{R}}$ with respect to the one-dimensional subspace $ \langle\pmb{z}_{2\nu+4}\rangle_\mathbb{R} $.
In particular, vectors $\pmb{f}_0$ and $\pmb{f}_2$ generate the same two-dimensional linear subspace
as the vectors $\pmb{z}_{2\nu+3}$  and $\pmb{z}_{2\nu+4}$. This is convenient for our purposes.
 Vectors $\pmb{f}_j, j =0,1,2$ are defined uniquely by these conditions. 
 
 We consider the image $\pmb{z}' = \frak{P}_2( \overline{\pmb{X}}_{2\nu+4})$ which is well-defined and for which we have
 $$
\pmb{z}' =  r' \pmb{f}_0+
 s'\pmb{f}_1+t' \pmb{f}_2 ,
 $$
where
 \begin{equation}\label{rst1w}
 r'\asymp \frac{q_{2\nu+4}^{1+\lambda \frak{g}_1}}{\Delta} ,\,\,\,\,\,\,
 s'= \frac{1}{\Delta},
 \,\,\,\,\,\,
 t' \asymp \frac{1}{\Delta}.
 \end{equation} 
 Note that  as $\overline{\pmb{X}}_{2\nu+4}\in S$, so does $ \pmb{z}' \in S$. From the last equality \eqref{g11} we see that 
   \begin{equation}\label{rst2w}
   r'\asymp q_{2\nu+4}^{\frak{g}_1}.
   \end{equation}
   The image $ \frak{P}_2( \overline{\mathcal{W}}_{2\nu+4 })$ is an ellipsoid. Its larger major axis is
   $\asymp r' $ and is larger than $ 6q_{2\nu+4}$. The section $ \{ x_0 = r'\}$ of the ellipsoid is a ball of radius 
   of order 
   $\asymp s' =\frac{1}{\Delta}\asymp q_{2\nu+3}^{-1+\lambda+\frac{\lambda}{\frak{g}_1}} $. From the right inequality from \eqref{prw} we see that $\xi_{2\nu+1} \le q_{2\nu+2}^{-\lambda} = o(s)$.
   So  the radius  of the section $ \{x_0 = r' \} $ of the ellipsoid is $ \ge 6\xi_{2\nu+1}$. By convexity
   $$ 
   \frak{P}_2( \overline{\mathcal{W}}_{2\nu+4 })\supset \pmb{z}' + \mathcal{C}
   ,$$
   and the lemma is proven as  $\mathcal{C}$ contains a fundamental domain of the lattice $\Gamma$.
   \end{proof}

From Lemma \ref{Lemma5} we see that  $ \frak{P}_2( \overline{\mathcal{W}}_{2\nu+4 })$ contains a shift of a fundamental domain of the lattice $\Gamma$ and the lattice $\Gamma_1$ is congruent to $\Gamma$. So $ \frak{P}_2( \overline{\mathcal{W}}_{2\nu+4 })$
 contains an integer point. This is just the point $\pmb{z}_{2\nu+5} = (q_{2\nu+5} , a_{1,2\nu+5}, a_{2,2\nu+5}a_{3,2\nu+5})$ which we are constructing. 
 
 As we defined $\pmb{z}_{2\nu+5}$, we may consider the angular neighborhood $\mathcal{U}_{2\nu+5}$.
 From the  inequalities on  parameters \eqref{rst1w} and triangle inequality we see that 
 $\mathcal{U}_{2\nu+5}\subset \mathcal{W}_{2\nu+4}$. Moreover, it is clear from the definition of the angular neighborhood
  $\mathcal{U}_{2\nu+5}$ that for all $\pmb{\eta} \in \mathcal{U}_{2\nu+5}$ vectors $\pmb{z}_{2\nu+4}$ and $\pmb{z}_{2\nu+5}$ are successive best approximation vectors to $\pmb{\eta}$.

 \vskip+0.3cm
  Now we will verify the rest of conditions ({\bf i}) - ({\bf vii}).
   
 \vskip+0.3cm
 
It is clear that the triple $\pmb{z}_{2\nu+3}, \pmb{z}_{2\nu+4}, \pmb{z}_{2\nu+5}$ is primitive and both quadruples 
$$\pmb{z}_{2\nu+1}, \pmb{z}_{2\nu+3}, \pmb{z}_{2\nu+4}, \pmb{z}_{2\nu+5}\textrm{ and } \pmb{z}_{2\nu+2}, \pmb{z}_{2\nu+3}, \pmb{z}_{2\nu+4}, \pmb{z}_{2\nu+5}$$ 
 form a basis of $\mathbb{Z}^4$.
 
 As for the value of $r'$ we have inequality \eqref{rst2w} and $\pmb{z}_{2\nu+5}\in \pmb{z}'+ \mathcal{C}$, we see that   
 $$
 q_{2\nu+5}\asymp q_{2\nu+4}^{\frak{g}_1},
 $$
 and we get (second) inequality from the condition ({\bf iv}) for the next step.
  
  So we have constructed two vectors $\pmb{z}_{2\nu+4}, \pmb{z}_{2\nu+5}$ which determine the next collection
  $$
   \pmb{z}_{2\nu+1}, \pmb{z}_{2\nu+2}, \pmb{z}_{2\nu+3},
  \pmb{z}_{2\nu+4}, \pmb{z}_{2\nu+5},
  $$
  and verified  conditions ({\bf i}) - ({\bf v}).     
  
  Geometric  condition ({\bf vi})  for the angle between vectors
  $\pmb{Z}_{2\nu+3}-\pmb{Z}_{2\nu+5}$ and 
  $\pmb{Z}_{2\nu+4}-\pmb{Z}_{2\nu+5}$ 
  is valid because of Lemma \ref{Lemma3} as we take $ \pmb{Z}_{2\nu+5} \in \mathcal{W}_{2\nu+4}$.
  
  As for ({\bf vii}), it immediately follows from ({\bf vi}) by means Lemma \ref{Lemma6}
  which we introduced in Subsection \ref{three}.
 Indeed, let us finish with ({\bf vii}).
  Indeed, let
  $\Delta_*$  be the fundamental volume 
  of the lattice
  $$
  \Gamma_* = 
   \langle\pmb{z}_{2\nu+3},\pmb{z}_{2\nu+4},\pmb{z}_{2\nu+5}
\rangle_\mathbb{Z}  = \mathbb{Z}^4 \cap  \langle\pmb{z}_{2\nu+3},\pmb{z}_{2\nu+4},\pmb{z}_{2\nu+5}
\rangle_\mathbb{R}
  .
 $$
 We apply Lemma \ref{Lemma6} to the vectors $\frak{z}_j =\pmb{z}_{2\nu+2+j}$ and take into account   ({\bf vi})  together with the inequalities 
 $$
  |\frak{Z}_j- \frak{Z}_3| =
  |\pmb{Z}_{2\nu+2+j} - \pmb{Z}_{2\nu+5}| \asymp \frac{ \zeta_{2\nu+2+j}}{q_{2\nu+2+j}},\,\,\,\, j = 1,2
  $$
  to obtain
  $$
  \Delta_* \gg q_{2\nu+2}q_{2\nu+4}q_{2\nu+5}  \cdot  \frac{ \zeta_{2\nu+3}}{q_{2\nu+3}}  \cdot \frac{ \zeta_{2\nu+4}}{q_{2\nu+4}} = \zeta_{2\nu+3}\zeta_{2\nu+4}q_{2\nu+5}.
  $$
  The upper bound 
   $$
  \Delta_* \ll \zeta_{2\nu+3}\zeta_{2\nu+4} q_{2\nu+5}
  $$
  is standard.
  
  We have checked all the conditions of the next step of induction. Thus, Theorem \ref{Thm2} in the case $k=1$ is proven.$\Box$
   
\section{Comments on the case $ k>1$}\label{9}
 
In the case $k=1$, we constructed vector $\pmb{\alpha}$ where best approximations form a periodic sequence of patterns $AB$, where growth of denominators $q_\nu$
 is determined by conditions ({\bf iv}) and ({\bf v}). That is, on the even step  ({\bf Stage 1}) the exponent of growth was equal to $\frak{g}_1$ and on the odd step  ({\bf Stage 2}) it is $\frak{g}_{1,0}$ defined in Lemma \ref{Lemma1} and \eqref{g11} respectively. The possibility of \emph{gluing}  pattern $A$ and pattern $B$ is enabled by equalities in  \eqref{g11} which comes from the equation \eqref{xRoy}. This alternance is enclosed in the sequence $u_j$ defining the sets $\mathcal{U}_j$, which is periodic with period of length $2$.
  To carry out the construction for arbitrary $k>1$ we simply need to construct periodic patterns 
   $ \underbrace{A\ldots A}_{k\, \text{times}}B$ 
   corresponding to triples \eqref{triples} with $ \nu_{j-1} =\nu_j +1,  1\le j \le k$. For this, first of all we should   define inductively
\begin{equation}\label{anyg}
\frak{g}_{k,0} = \frac{1}{\theta\frak{g}_{k} -1},\,\,\,\,\,
\frak{g}_{k,j} = f_j(\frak{g}_{k},\lambda),\, j = 1,\ldots, k-1,
\end{equation}
where functions $f_j(\cdot,\cdot)$ are defined in \eqref{fu2}.
Remark 4 will be of importance.
Calculations show that 
\begin{equation}\label{ssdd}
\max_{0\le j \le k-1} \frak{g}_{k,j} \le \frak{g}_{k}.
\end{equation}
Equalities \eqref{india} and \eqref{anyg} give us
\begin{equation}\label{glue}
\frak{g}_k= \frac{\frak{g}_{k,0}+1}{\theta\frak{g}_{k,0}}
,\,\,\,\,\,
\frak{g}_{k,j-1}
 = \frac{\frak{g}_{k,j}-\lambda}{(1-\lambda) \frak{g}_{k,j }},\,\,\,\
 j = 1, \ldots, k-1.
 \end{equation}

   Then we consider the sets $\mathcal{U}_j$ defined by a periodic sequence $u_j$ with period of length $k+1$, defined by
    $$
 u_j =
 \begin{cases}\label{uj}
 1+ \lambda \frak{g}_{k,l},\,\,\,\,\,\,\, j \mod k+1 = l
  \cr
 1+ \lambda \frak{g}_k, \,\,\,\,\,\, j \mod k+1 =k
 \end{cases}
 $$
   where $\frak{g}_k$ is defined in Lemma \ref{Lemma1} while $\frak{g}_{k,l}$ for $0\le l \le k-1$ are defined in \eqref{anyg}. In the inductive construction, at step $j$ we apply ({\bf Stage 2}) if $ j \mod k+1 =k$ and ({\bf Stage 1}) otherwise. In this construction the last two elements of each triple
 $ \pmb{z}_{\nu_j-1}, \pmb{z}_{\nu_j}, \pmb{z}_{\nu_j+1}$ 
 are just the first two elements of the triple
  $ \pmb{z}_{\nu_{j-1}-1}, \pmb{z}_{\nu_{j-1}}, \pmb{z}_{\nu_{j-1}+1}$,
  that is $\nu_{j-1} = \nu_j +1$.
 The possibility of gluing together separate  patterns $A$ and $B$ into 
 $ \underbrace{A\ldots A}_{k\, \text{times}}B$  follows from \eqref{glue} because of the identities
 $$
 q_{\nu_j+1}^{(1-\lambda)\frak{g}_{k,j}-1} \asymp \frac{\Delta}{d_j}
, $$
where
$
\Delta
$ is the height of the three-dimensional rational subspace containing
all the triples $ \pmb{z}_{\nu_j-1}, \pmb{z}_{\nu_j}, \pmb{z}_{\nu_j+1},
0\le j \le k$ and
$
d_j \asymp
 q_{\nu_j+1}^{1-\lambda}
 $
 is the height of  the two-dimensional subspace  generated by vectors
 $ \pmb{z}_{\nu_j}, \pmb{z}_{\nu_j+1}$.

 The conclusion $ \frac{\omega(\pmb{\alpha})}{\hat{\omega}(\pmb{\alpha}) } =\frak{g}_k
$ follows from \eqref{ssdd}.\\

 {\bf Acknowledgment}
 The authors would like to thank the \emph{Centro Internazionale per la Ricerca Matematico} (CIRM) in Trento for the opportunity of a \emph{Research in Pairs} stay from May 19th to June 1st 2019. Many ideas of this paper appeared during this fruitful time.

Research is supported by joint FWF-Projekt I 5554 and RSF project 22-41-05001 "Diophantine approximations, arithmetic sequences \& analytic number theory" (https://rscf.ru/en/project/22-41-05001/).
% The authors want to thank the anonymous referee for the very careful reading of the manuscript, many valuable comments and essential  simplification of the exposition.

 \vskip+1cm
 
 authors:
 \vskip+0.3cm

 Antoine Marnat,
 
Univ Paris Est Creteil, Univ Gustave Eiffel, CNRS, LAMA UMR8050, F-94010 Creteil, France\\
\indent Institut f\"ur Diskrete Mathematik und Geometrie, TU Wien,  Wien, Austria
 
 \vskip+0.3cm
 
 Nikolay Moshchevitin,
 
Faculty of Computer Science, HSE University, Pokrovsky boulevard 11, Moscow, Russia 109028 
   
\end{document}